\documentclass[journal]{IEEEtran}
\usepackage{bm}
\usepackage{amsmath,amsfonts,amssymb}
\usepackage{algorithmic}
\usepackage{soul,color}
\usepackage{epsfig}
\usepackage{graphicx,float,wrapfig,subfigure}
\usepackage{cite}
\usepackage{epstopdf}
\usepackage{url}
\usepackage[colorlinks,linkcolor=black,anchorcolor=black,citecolor=black,urlcolor=black]{hyperref} 
\usepackage[normalem]{ulem}
\usepackage{makecell}

\usepackage[ruled,linesnumbered]{algorithm2e}

\newtheorem{colloary}{Corollary}
\newtheorem{theorem}{Theorem}
\newtheorem{proposition}{Proposition}
\newtheorem{lemma}{Lemma}

\newtheorem{remark}{Remark}

\usepackage{threeparttable}

\usepackage{color}
\makeatletter %
\let\myorg@bibitem\bibitem
\def\bibitem#1#2\par{%
	\@ifundefined{bibitem@#1}{%
		\myorg@bibitem{#1}#2\par
	}{%
		\begingroup
		\color{\csname bibitem@#1\endcsname}%
		\myorg@bibitem{#1}#2\par
		\endgroup
	}%
}

\makeatother %

\allowdisplaybreaks

\usepackage{xcolor}
\newcommand{\cg}[1]{\textcolor{black}{#1}}
\begin{document}

\title{
Neural Risk Limiting Dispatch in Power Networks: \\
Formulation and Generalization Guarantees
	}

\author{
Ge~Chen,~\IEEEmembership{Member,~IEEE,}
and~Junjie~Qin,~\IEEEmembership{Member,~IEEE}
\thanks{
This paper is supported in part by the U.S.  National Science Foundation under Grant No. ECCS-2339803 and under Grant No. ECCS-2434502. Earlier versions of this paper benefited from helpful discussion with  Kameshwar Poolla and Pravin Varaiya. (Corresponding author: \textit{Junjie Qin})

G. Chen is with the School of Engineering at Great Bay University, Dongguan, China, and also with the Elmore Family School of Electrical and Computer Engineering at Purdue University, West Lafayette, Indiana, USA. J. Qin is with the Elmore Family School of Electrical and Computer Engineering at Purdue University, West Lafayette, Indiana, USA. (Emails: chen4911@purdue.edu, jq@purdue.edu)}
	
}

\maketitle

\begin{abstract}
Risk limiting dispatch (RLD) has been proposed as an approach that effectively trades off economic costs with operational risks for power dispatch under uncertainty. However, how to solve the RLD problem with provably near-optimal performance still remains an open problem. This paper presents a learning-based solution to this challenge. We first design a data-driven formulation for the RLD problem, which aims to construct a decision rule that directly maps day-ahead observable information to cost-effective dispatch decisions for the future delivery interval. Unlike most existing works that follow a predict-then-optimize paradigm, this end-to-end rule bypasses the additional suboptimality introduced by separately handling prediction and optimization. We then propose neural RLD, a novel solution method to the data-driven formulation. This method leverages an L2-regularized neural network (NN) to learn the decision rule, thereby transforming the data-driven formulation into a training task that can be efficiently completed by stochastic gradient descent. A theoretical performance guarantee is further established to characterize the suboptimality of our method, which implies that its suboptimality approaches zero with high probability as more samples are utilized. Simulation tests across various systems demonstrate our method's superior performance in convergence, suboptimality, and computational efficiency compared with benchmarks.   
\end{abstract}

\begin{IEEEkeywords}
Risk limiting dispatch, renewable integration, deep learning, stochastic programming,  generalization performance guarantee. 
\end{IEEEkeywords}

\section{Introduction} \label{sec_intro}

The increasing adoption of uncertain renewables necessitates methodologies for power dispatch under uncertainty to achieve safe and cost-effective power system operations \cite{javed2023impact}. These methodologies typically fall into two major categories: i) robust optimization that safeguards the system performance for the worst case \cite{zhao2022sustainable} and ii) stochastic programming and stochastic control formulations that optimize the expected performance  \cite{zakaria2020uncertainty}.

Within the second category, Risk Limiting Dispatch (RLD) distinguishes itself as a class of methods that effectively trade off the economic cost with operation risks utilizing the most updated information about the uncertainties \cite{5618534, rajagopal2013risk}. Varaiya et al. \cite{5618534} initially proposed RLD within a multi-stage stochastic control formulation to incorporate recourse decisions in multi-settlement markets. This approach was later generalized to explicitly leveraging forecast updates (i.e., the progressively improving forecasts for loads and renewable generation as we approach the delivery time) \cite{rajagopal2013risk}, generator ramping constraints \cite{qin2013risk}, and energy storage \cite{6580485}. These early studies mostly relied on deriving analytical solutions for the underlying stochastic control problems, circumventing the \emph{curse of dimensionality} associated with continuous state/action/disturbance variables in the original RLD formulation. However, when power network constraints are considered (e.g., as in the network RLD problem), multi-stage stochastic control formulation for RLD turns out to be intractable due to the high-dimensional continuous state space and the lack of analytical solutions. In fact, even the two-stage special case of the network RLD problem is challenging, for which {Zhang et al. \cite{6818365} proposed a heuristic algorithm and established a performance guarantee to characterize its suboptimality under a ``small-$\sigma$ assumption", where $\sigma$ represents the standard deviation of forecast errors. This assumption implies that $\sigma$ is much smaller than the expected net demand. Under this condition, forecast errors have minimal impact on the network's operational characteristics, and the dispatch decisions made under both deterministic and stochastic settings remain similar.} 
It is noteworthy that the formulation and proposed algorithm in \cite{6818365} do not rely on sampling from the probability distribution of the disturbances, which can be computationally costly when the disturbances are high-dimensional (e.g., considering correlated spatially dispersed wind and solar outputs). This makes their proposed algorithm easy to be incorporated in existing power system dispatch software (e.g., by modifying the reserve margin based on RLD solutions).  {However, this assumption may not hold in practice, particularly in power systems with high renewable penetration, where forecast errors can be significant. In such cases, the performance guarantee established under the small-$\sigma$ assumption may no longer apply, especially when renewable generation results in negative net demand, as the authors of \cite{6818365} acknowledge.} Hence, it remains an \emph{open problem to develop provably near-optimal methods that solve the RLD problem with network constraints and practical uncertainty levels} (i.e., without the small-$\sigma$ assumption).

Another challenge associated with RLD and power dispatch under uncertainty in general is to obtain accurate characterizations for the uncertainties. Stochastic control formulations of RLD require the joint probability distribution of all uncertain parameters involved, which may be difficult to estimate from data when a large number of distributed wind and solar generators are considered \cite{8618917}. Probabilistic forecasting methods based on machine learning \cite{9599541, 8982039} could be leveraged to estimate these joint distributions. Other formulations have been proposed to address power dispatch under uncertainty, including scenario-based methods which represent the high-dimensional joint distributions using discrete scenarios \cite{8999581, 8651522}. However, theoretical understanding of the performance of these formulations for characterizing uncertainty is lacking. Even when guarantees are available, it remains unclear how errors in quantifying the uncertainty may impact subsequent dispatch tasks derived from the uncertainty model.

Even if the difficulties associated with characterizing the uncertainties, i.e., the \emph{prediction} (albeit probabilistic) step, and identifying the optimal decisions given the uncertainties, i.e., the \emph{optimization} step, can be resolved, fundamental pitfalls around the two-step \emph{predict-then-optimize} paradigm make it challenging to achieve direct data-driven, end-to-end optimal or provably near optimal RLD decisions. Indeed, it is common for the prediction step to use standard statistical loss functions, e.g., mean squared error (MSE), to quantify the accuracy ignoring the actual objective functions in the subsequent optimization step, and the optimization step assumes the uncertainty model obtained from the prediction step is precise.
However, when the errors in the prediction step are accounted for, a ``more accurate'' model for uncertainty as measured by MSE may not always lead to a better decision in the optimization step \cite{elmachtoub2022smart, 9755891}. Because of this ``mismatch" in the objectives of the two steps, this paradigm may introduce additional suboptimality, resulting in undesirable performance.

\subsection{Contributions and paper organization}
In this paper, we propose \emph{neural RLD}, a formulation and method that merges the prediction and optimization steps, which directly use the data to train a neural network (NN) outputting dispatch decisions. This enables us to circumvent the challenges associated with the individual steps in the overall pitfalls of the predict-then-optimize paradigm, while obtaining end-to-end generalization bounds for the quality of the resulting dispatch decisions. This paper contributes to the literature in the following ways:
\begin{enumerate}
\item We design a data-driven formulation for the network RLD problem. This formulation aims to learn a decision rule from historical data. This decision rule is end-to-end, meaning it directly maps day-ahead observable information to cost-effective dispatch decisions for future delivery intervals, thereby bypassing the ``mismatch" issue of the predict-then-optimize paradigm.
\item We propose the neural RLD method to solve the data-driven formulation above. This method utilizes an L2-regularized NN to learn the decision rule, effectively converting the data-driven formulation into a NN training task with a specialized training loss. We then explicitly derive the gradient of this specialized loss, enabling efficient completion of this training task by using stochastic gradient descent (SGD).
\item {We establish a performance guarantee that theoretically characterizes the generalization error of the proposed method. This guarantee demonstrates that our method is asymptotically optimal if the NN is well-trained, implying that its optimality gap approaches zero with high probability as the sample size increases. To the best of our knowledge, this is the first performance guarantee developed for an end-to-end RLD approach.}
\end{enumerate}

The remaining parts are organized as follows. Section \ref{sec_literature} reviews related literature.
Section \ref{sec_modeling} describes the problem formulation of the network RLD. Section \ref{sec_solution} introduces the proposed neural RLD. Section \ref{sec_case} demonstrates simulation results, and Section \ref{sec_conclusion} concludes this paper.

\subsection{Other related literature} \label{sec_literature}

\subsubsection{Learning-based stochastic dispatch}

The application of machine learning to stochastic dispatch presents both innovative solutions and challenges. Initially, machine learning methods were usually deployed in one of two steps in the predict-then-optimize paradigm: enhancing prediction accuracy in the prediction step \cite{zhang2021review}, or improving computational efficiency in the optimization step \cite{10058008}. However, the aforementioned mismatch issue may bring additional suboptimality. {To address this limitation, end-to-end frameworks have emerged as a promising alternative. For example, imitation learning trains NNs to directly map observed features to dispatch decisions, thereby avoiding the mismatch issue. However, this framework typically uses mean squared error (MSE) between predicted and actual optimal decisions as the training loss, which may not effectively capture the suboptimality of the predicted decisions \cite{elmachtoub2022smart}. Some studies have instead used the original objective of the power dispatch problem as the training loss, which can directly measure suboptimality \cite{10700765,10256159}. Nevertheless, most of these studies are confined to deterministic power dispatch tasks. A few studies have developed learning-based end-to-end methods for stochastic problems \cite{NIPS2017_3fc2c60b,pmlr_v202_rychener23a}, relying on limited historical data to approximate the true distribution of uncertain operation conditions in future delivery intervals. This discrepancy between the empirical distribution derived from data and the true population distribution introduces a generalization error that can induce suboptimality. Yet, none of these studies has provided a theoretical performance guarantee to quantify this generalization error and assess its potential impact on suboptimality, thereby limiting their practical applicability.}

%

\subsubsection{Sample average approximation}
Sample average approximation (SAA) is a method \cite{doi:10.1137/S1052623499363220} that uses a finite number of discrete scenarios to approximate uncertainties, which can effectively transform RLD into a deterministic form. 
Some error bounds have also been proposed to explain its performance on characterizing uncertainty. For example, Shapiro et al. \cite{shapiro2021lectures} have shown that the SAA solution uniformly converges to the optimal one as the number of scenarios increases indefinitely. Additionally, Kim et al. \cite{kim2015guide} established a link between the number of scenarios used and the worst-case optimality gap of the SAA. However, the SAA still suffers from the aforementioned ``mismatch" issue, as it is typically incorporated in the optimization step of the predict-then-optimize paradigm.

\subsubsection{Direct learning}

In conference paper \cite{8618917}, we initially proposed an end-to-end learning method to predict RLD decisions directly from historical data. This method offered a generalization guarantee to characterize its suboptimality. 
However, the initial study was limited to scalar decision variables, which significantly restricted its applicability and practicality in scenarios involving power network constraints. Meanwhile, it focused on theoretical analysis without numerical case studies. The current paper significantly extends the scope of \cite{8618917} by concentrating on the network RLD with vector variables, providing a new generalization performance guarantee for this networked case, and conducting a series of case studies to validate the practical applicability of our method.

\section{Problem Formulation} \label{sec_modeling}

\begin{figure}[h]

    \vspace{-4mm}
	\centering
	{\includegraphics[width=0.95\columnwidth]{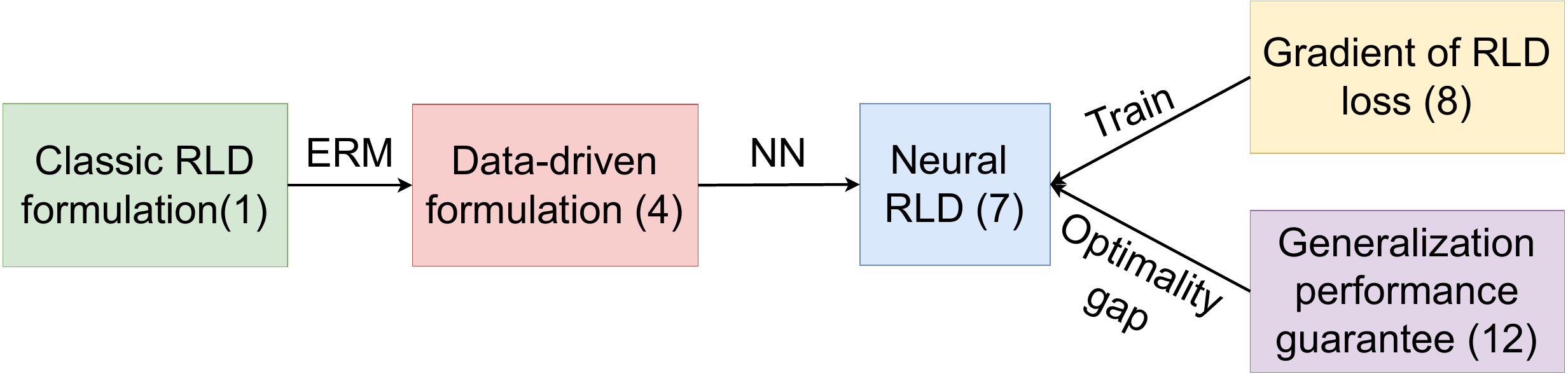}}
	\vspace{-2mm}
 	\caption{Steps for deriving the proposed neural RLD formulation.}
	\label{fig_summary}
	\vspace{-2mm}
\end{figure}

{This paper proposes a neural RLD approach, where a NN is trained to learn a decision rule that maps observable day-ahead information to cost-effective generation decisions. The detailed procedure for deriving the neural RLD formulation is summarized in Fig. \ref{fig_summary}. We start by presenting the classic network RLD \eqref{eqn_RLD} and the data-driven formulation \eqref{eqn_ERM} derived using empirical risk minimization (ERM). Then, the data-driven formulation is translated into a NN training task, resulting in the proposed neural RLD \eqref{eqn_h_hat}. The gradient of the training loss is explicitly derived using the envelope theorem, enabling efficient NN training. Furthermore, a generalization performance guarantee is established using statistical learning theory, which theoretically characterizes the suboptimality of our neural RLD.} 

This section introduces the classic and data-driven RLD formulations. The remaining methodology will be detailed in Section \ref{sec_solution}.

\subsection{Classic network risk limiting dispatch}

\def\Nbus{{N_\mathrm{B}}}
\def\Nline{{N_\mathrm{L}}}

As in \cite{6818365}, we consider a two-stage formulation of network RLD, where  the system operator (SO) makes dispatch decisions for a future delivery interval sequentially in the \emph{day-ahead} and \emph{real-time} stages. 
Detailed discussion on how this formulation fits with the two-settlement electricity markets and the practical dispatch procedures can be found in \cite{5618534, 6818365}. 
 
Consider a power network with $\Nbus$ buses. In the day-ahead stage, the SO dispatches generators based on the most updated information regarding the net demand (i.e., load and renewable generation) over the network. Denote this dispatch decision by $\mathbf u \in \mathbb R^\Nbus$, where $u_n$ is the generation at bus $n$. As the actual net demand $\mathbf d\in \mathbb R^\Nbus$ is not available at day ahead, the classical RLD formulation assumes that the probability distribution of $\mathbf d$ can be estimated. This probability distribution can depend on various auxiliary factors (e.g., temperature and wind speed information from the day-ahead weather forecast), and thus is denoted by $\mathcal D(\mathbf x)$, where $\mathbf x \in \mathbb R^p$ is the $p\times 1$ vector of auxiliary information available at day ahead. In the real-time stage, the SO observes the realization of net demand $\mathbf d$ and can adjust the day-ahead dispatch by either dispatching extra generation or reducing the day-ahead dispatched amount, so the adjusted total generation profile balances the grid. We denote the real-time adjustment by  $\mathbf g \in \mathbb R^\Nbus$. Mathematically, the two-stage network RLD problem takes the form of 
 \begin{align} 
	&\min_{\mathbf u} \quad \bm \alpha^\intercal \mathbf u + \mathbb{E}_{\mathbf d \sim \mathcal{D}(\mathbf x)}[Q(\mathbf u, \mathbf d)], \label{eqn_RLD}
\end{align}
where $\bm \alpha \in \mathbb R^\Nbus$ is the vector of linear generation cost coefficients for the day-ahead stage, and $Q(\mathbf u, \mathbf d)$ is the \emph{optimal} second-stage cost given the day-ahead dispatch $\mathbf u$ and the realized net demand $\mathbf d$, defined as the optimal value of the real-time dispatch problem. In particular, as in \cite{6818365}, we use the following optimal DC (i.e., linearized AC) power flow formulation for the real-time dispatch: 
\begin{subequations} \label{eqn_RLD_2nd}
\begin{align} 
	Q(\mathbf u, \mathbf d) := \min_{\mathbf g, \bm \theta} \quad &\bm \beta^\intercal (\mathbf g)_+, \label{eqn_obj_2nd}  \\
		 \mbox{s.t.} \quad & \mathbf u + \mathbf g - \mathbf d = \mathbf B \bm \theta\qquad\, : \bm \lambda, \label{eqn_DCOPF} \\
  & -\mathbf f^\mathrm{max} \leq \mathbf F \bm \theta \leq \mathbf f^\mathrm{max}:\bm \nu_{-},\bm \nu_{+}. \label{eqn_capacity}
\end{align}	
\end{subequations}
{Objective \eqref{eqn_obj_2nd} models the cost of generation adjustments or voluntary load shedding (like demand response) in real time, where $(\mathbf g)_+ := \max\{\mathbf g, \mathbf 0\}$ and $\bm \beta \in \mathbb{R}^\Nbus$ represents the per unit cost.
Specifically, a positive component $g_n$ in $\mathbf g$ corresponds to two possible scenarios: i) Generator $n$ increases its output from $u_n$ to $u_n + g_n$, incurring an additional generation cost of $\beta_n g_n$, or ii) power demand at bus $n$ is curtailed from $d_n$ to $d_n - g_n$, incurring a penalty of $\beta_n g_n$. Thus, the term $\bm \beta^\intercal (\mathbf g)_+$ captures the extra cost associated with these real-time actions.} Eq. (\ref{eqn_DCOPF}) represents the DC power flow equation, where matrix $\mathbf B \in \mathbb R^{\Nbus\times \Nbus}$ maps voltage phase angles $\bm \theta \in \mathbb R^\Nbus$ into nodal active power injections. Eq. (\ref{eqn_capacity}) denotes the branch capacity constraints, where $\mathbf F \in \mathbb R^{\Nline \times \Nbus}$ maps $\bm \theta$ to branch flows, $\mathbf f^\mathrm{max}$ is the corresponding maximum allowable values, and $\Nline$ is the branch number. Vectors $\bm \lambda \in \mathbb R^\Nbus$, $\bm \nu_{-} \in \mathbb R^\Nline$  and $\bm \nu_{+} \in \mathbb R^\Nline$ are the dual variables of (\ref{eqn_DCOPF})-(\ref{eqn_capacity}), respectively. 
{We use the DC power flow formulation instead of the AC formulation because, in practical wholesale markets, day-ahead generation decisions in transmission networks are typically made based on the DC model \cite{1266587}. Moreover, the accuracy of the DC formulation is generally acceptable in transmission network operations, as line reactances are relatively small compared to resistances \cite{li2025frequency}.}

\begin{remark}[VOLL interpretation]
Alternatively, the objective \eqref{eqn_obj_2nd} can be interpreted as the \emph{value of lost load} (VOLL), in which case $\bm \beta$ is the per unit cost that users are willing to pay to avoid a loss of load event (involuntary load shedding). As a result, $\mathbb E_{\mathbf d}[Q (\mathbf u, \mathbf d)]$ models the \emph{risk} associated with potential \emph{loss of load} events\footnote{VOLL is among one of the widest used power system risk metrics \cite{kariuki1996evaluation}. A power system operator that utilizes the objective in \eqref{eqn_RLD} to make its decision is \emph{risk-averse} in the sense of \cite[Definition 6.C.1]{mas1995microeconomic} as it is easy to show that $ Q(\mathbf u, \mathbb E\mathbf d) \le  \mathbb E Q(\mathbf u, \mathbf d) $ for any $\mathbf u$.}, and \eqref{eqn_RLD} trades off day-ahead generation cost with real-time operational risk. Note the values of $\bm \beta$ here are much higher compared to the voluntary load shedding scenario. Furthermore, the power system operator can reduce the risk of loss of load event by choosing larger values for $\bm \beta$.
\end{remark}


 The network RLD problem is typically tackled by a two-step predict-then-optimize paradigm, i.e., first estimating the distribution $\mathcal{D}$ via probabilistic forecasting, and then solving the stochastic problem \eqref{eqn_RLD}. However, as mentioned in Section \ref{sec_intro}, separately considering the prediction and optimization leads to challenges for each step and loss of optimality overall.

\subsection{Data-driven limiting dispatch}

We next introduce a new formulation of two-stage network RLD that merges the prediction and optimization steps in the predict-then-optimize paradigm. In other words, we will leverage the historical data that are commonly used in the prediction step to directly learn a decision rule for the future delivery interval. We refer to this formulation as \emph{data-driven network RLD}.

The historical dataset includes day-ahead auxiliary information $\mathbf x^{(m)}$ and the realized net demand $\mathbf d^{(m)}$ for a collection of past delivery intervals $m \in \mathcal M$ with $M:=|\mathcal M|$.  As this dataset $\mathcal S:= \{(\mathbf x^{(m)}, \mathbf d^{(m)}), m \in \mathcal M\}$ will be used to train a decision rule for future delivery intervals, we refer to $\mathcal S$ as the \emph{training set}. The training set  is usually available to the SO well in advance of the delivery interval for which the SO needs to make dispatch decisions about. At the day-ahead stage for the delivery interval of concern, the SO will also have access to the auxiliary information for the particular delivery interval, denoted by $\mathbf x'$. In the prediction step of the conventional predict-then-optimize paradigm, if a deterministic forecast is used in place of a probabilistic one for the future net demand, $\mathbf x'$ will be used together with a machine learning model trained on $\mathcal S$ for demand forecast. In our setting, we wish to apply the decision rule trained from $\mathcal S$ at $\mathbf x'$ to predict the RLD decisions. We refer to $\mathbf x'$ as the \emph{test input}.

Ideally, we wish to find the optimal decision rule $\mathbf h^\star: \mathbb R^p \mapsto \mathbb R^\Nbus$, that maps the day-ahead auxiliary information $\mathbf x$ to the optimal dispatch decision $\mathbf u^\star$.
Adapting terminologies from the statistical learning theory, we will also refer to such a decision rule as a \emph{hypothesis}. Mathematically, supposing that there is an underlying true joint distribution of $(\mathbf x, \mathbf d)$, denoted by $\mathcal Z$, we can write the following problem to identify the optimal hypothesis as
\begin{align}
\mathbf h^\star = \mathop{\arg\min}_{\mathbf h \in \mathcal{H}} \ \mathbb{E}_{(\mathbf{x}, \mathbf{d}) \sim \mathcal{Z}} \left[\bm \alpha^\intercal \mathbf h(\mathbf x) + Q(\mathbf h(\mathbf x),\mathbf d)\right], \label{eqn_optimal}
\end{align}
where $\mathcal{H}$ is a hypothesis class containing all functions that map $\mathbf{x}$  to $\mathbf{u}$. 
However, directly evaluating the expectation in \eqref{eqn_optimal} is challenging as the distribution $\mathcal Z$ is usually not available. Instead, we approximate the expectation with its empirical mean using the training set $\mathcal S$ as done in empirical risk minimization (ERM) \cite{montanari2022universality}: 
\begin{align}
\widehat{\mathbf h} = \mathop{\arg\min}_{\mathbf h \in \mathcal{H}} \frac{1}{M} \sum_{m \in \mathcal{M}} \left(\bm \alpha^\intercal \mathbf h(\mathbf x^{(m)}) + Q(\mathbf h(\mathbf x^{(m)}),\mathbf d^{(m)})\right).  \label{eqn_ERM}
\end{align}
We refer to \eqref{eqn_ERM} as the data-driven RLD problem. Upon solving \eqref{eqn_ERM}, and at the day-ahead stage for the future delivery time of concern, we can directly obtain the dispatch decision via $\mathbf u' = \widehat{\mathbf h} (\mathbf x')$.  

This data-driven formulation bypasses the ``mismatch'' issue discussed previously by merging the prediction and optimization steps. Moreover, it does not need data of the optimal decision $\mathbf u^\star$ for historical delivery interval $m\in \mathcal M$ as the training labels. However, it poses challenges for two primary reasons. First, given that $\mathbf h \in \mathcal{H}$ can represent any  $\mathbb R^\Nbus$-valued function, Eq. \eqref{eqn_ERM} is computationally intractable due to its infinite-dimensional nature. Second, we only have a finite number of historical data points in the training set. These data points may fail to accurately reflect the true distribution of $\mathcal{Z}$, so the ERM hypothesis $\widehat{\mathbf h}$ may perform much worse than $\mathbf{h}^\star$.

\section{Neural Risk Limiting Dispatch} \label{sec_solution}

Motivated by the above challenges, we propose the neural RLD. In this approach, the hypothesis class $\mathcal{H}$ is represented by L2 regularized NNs with a finite number of parameters. Consequently, problem \eqref{eqn_ERM} is redefined as a finite-dimensional NN training task with a specialized training loss. To support this approach, we provide a generalization performance guarantee based on statistical learning tools including uniform convergence and Rademacher complexity, which theoretically characterizes the suboptimality of our method.\footnote{It is impossible to obtain performance guarantees here with probability one, as there is always a non-zero probability (as long as the training set is finite) to obtain a training set for which the empirical mean in \eqref{eqn_ERM} is not a good representation of the population mean in \eqref{eqn_optimal}.}  This guarantee demonstrates that as more historical data points are utilized, the performance of our neural RLD approaches that of the optimal hypothesis in the same hypothesis class. In this section, we first describe the overall neural RLD procedure. Following this, we delve into the training process of the neural RLD. Finally, we establish the generalization performance guarantee.

\subsection{Neural network for data-driven RLD}\label{sec_E2E}
{The mapping from $\mathbf x$ to $\mathbf u^\star$ is nonlinear in general because: i) the net demand $\mathbf d$ may not have a linear relationship with the day-ahead observable information $\mathbf x$, and ii) even if $\mathbf d$ is linearly dependent on $\mathbf{x}$, the mapping from $\mathbf x$ to $\mathbf u^\star$ is not necessarily linear because it involves an optimization problem.} 
\cg{Hence, we select NNs as the hypothesis class due to their superior approximation capability \cite{NEURIPS2020_2000f632}. To mitigate overfitting, L2 regularization is employed. Compared to L1 regularization, L2 regularization promotes smoother and more stable optimization by penalizing large weights proportionally to their square, which also facilitates faster convergence during training.} 
By using $k \in \mathcal{K}$ and $j \in \mathcal{J}^{k}$ with $K:=|\mathcal K|$ and $J^k:=|\mathcal{J}^k|$ to index the hidden layers and the neurons in each, this hypothesis class can be expressed as:\footnote{\cg{Unlike many studies that incorporate L2-regularization by adding a penalty term to the training loss, we employ the constraint \(\|\mathbf{w}^{k}_{j}\|_2 \leq W^{\max}\) to simplify the derivation of the proposed performance guarantee. Furthermore, as noted in \cite[Section 3.4.1]{hastie2017elements}, regularization through penalty terms and constraints is theoretically equivalent.
}} 
\begin{align}
\mathcal{H} = \left\{
\mathbf h(\mathbf x) \left|
\begin{aligned}
&\mathbf h(\mathbf x) = \mathbf W^K\phi \left(\cdots \mathbf W^{2} \phi\left(\mathbf W^{1} \mathbf x\right)\right), \\
&\bm \Vert \mathbf w^{k}_{j} \Vert_2 \leq W^{\max}, \forall j \in \mathcal{J}^{k}, \forall k \in \mathcal{K},
\end{aligned} \right. \label{eqn_class}
\right\},
\end{align}
where $\mathbf h: \mathbb R^p \mapsto \mathbb R^\Nbus$ represents a NN in this hypothesis class; $\mathbf W^{k} \in \mathbb{R}^{J^{k+1} \times J^{k}}$ is the weights of the $k$-th layer and $\mathbf{w}_{j}^k \in \mathbb{R}^{J^k}$ is the transpose of its $j$-th row. Symbol $\phi(\cdot)$ denotes the activation function. This paper chooses ReLU as the nonlinear activation function, i.e., $\phi(\mathbf x)=\max\{\mathbf x, \bm 0\}$. Hence, this NN can represent nonlinear decision rules. 
It operates on each individual element of the input, ensuring that the dimensions of both input and output remain identical. The second constraint in \eqref{eqn_class} regularizes the class of NNs by limiting the L2 norm of the weights, where $W^{\max}$ is the maximum norm of weights allowed. If we further define a new loss function $\ell$, also referred as to the RLD loss, as
\begin{align}
\ell(\mathbf h(\mathbf x), \mathbf d) = \bm \alpha^\intercal \mathbf h(\mathbf x) + Q(\mathbf h(\mathbf x),\mathbf d), \label{eqn_loss}
\end{align}
then problem \eqref{eqn_ERM} can be interpreted as a NN training task. This task uses $\mathcal{S} = \{(\mathbf{x}^{(m)}, \mathbf{d}^{(m)}), \forall m \in \mathcal{M}\}$ as the training set, and leverages $\ell$ as the training loss, as follows:
\begin{align}
\widehat{\mathbf h} = \mathop{\arg\min}\limits_{h \in \mathcal{H}} \frac{1}{M}\sum_{m \in \mathcal{M}} \ell(\mathbf h(\mathbf x^{(m)}), \mathbf d^{(m)}). \label{eqn_h_hat}
\end{align}

\subsection{NN training} \label{sec_gradient}
The training task \eqref{eqn_h_hat} can be effectively completed using the stochastic gradient descent (SGD) algorithm if the (sub-) gradient of the RLD loss \eqref{eqn_loss} with respect to the weights exists and can be effectively evaluated. Based on the chain rule, this gradient is expressed as:
\begin{align}
 \nabla_{\mathbf{w}_{j}^k} \ell = \left( \bm{\alpha} + \nabla_{\mathbf{u}} Q \right)^\intercal  \nabla_{\mathbf{w}_{j}^k} \mathbf{u}, \forall j \in \mathcal{J}^{k+1}, \forall k \in \mathcal{K},  \label{eqn_gradient_full}
\end{align}
where $\mathbf u \in \mathbb{R}^{\Nbus}$ is the output of the hypothesis $\mathbf{h}$, and we have $\nabla_{\mathbf{w}_{j}^k} \ell \in \mathbb{R}^{J^k}$, $\nabla_{\mathbf{u}} Q \in \mathbb{R}^{\Nbus}$, and $\nabla_{\mathbf{w}_{j}^k} \mathbf{u} \in \mathbb{R}^{\Nbus \times J^k}$. The gradient $\nabla_{\mathbf{w}_{j}^k} \mathbf{u}$ can be easily obtained using automatic differentiation in modern deep learning frameworks like Pytorch. However, the gradient $\nabla_{\mathbf{u}} Q$ is hard to calculate because $Q(\mathbf u,\mathbf d)$ involves an optimization problem \eqref{eqn_RLD_2nd}. To address this issue, we leverage the strong duality to explicitly express $\nabla_{\mathbf{u}} Q$, as described in the following proposition.

\begin{proposition}[Expression of gradients] \label{prop_1} 

The gradient $\nabla_{\mathbf{u}} Q$ can be expressed as \cite{liang2022operation}:
\begin{align}
\nabla_{\mathbf{u}} Q = \bm \lambda^\star, \label{eqn_gradient}
\end{align}
where $\bm \lambda^\star$ is the optimal dual variable of constraint \eqref{eqn_DCOPF}.

\end{proposition}

\emph{Proof}: See Appendix \ref{app_3}.

Then, we can introduce the SGD to train the NN.
The whole training process is summarized in Algorithm \ref{algo_1}. Here, we use $N^\mathrm{ep}$ and $\mathcal{B}$ to denote the maximum epoch number and the sample index of a batch, respectively. 

\begin{algorithm} \label{algo_1}
  \SetKwInOut{Input}{Input}\SetKwInOut{Output}{Output}
  \SetAlgoLined
  \Input{Dataset $\mathcal{S} = \{(\mathbf{x}^{(m)}, \mathbf{d}^{(m)}), \forall m \in \mathcal{M}\}$, learning rate $\rho$, and batch size $B = |\mathcal{B}|$}
  \Output{The hypothesis ${\mathbf h}$}
  Initialize weights $\mathbf W^{k}, \forall k \in \mathcal{K}$\;
  \For{$epoch \gets 1$ to $N^\mathrm{ep}$}{
    \ForEach{batch $\{(\mathbf{x}^{(m)}, \mathbf{d}^{(m)}), \forall{m} \in \mathcal{B}\}$}{${\mathbf u}^{(m)} \gets \mathbf h({\mathbf{x}}^{(m)}), \forall m \in \mathcal{B}$\;
    Solve the dual of \eqref{eqn_RLD_2nd} to get $(\bm \lambda^\star)^{(m)}, \forall{m} \in \mathcal{B}$\;
    $\mathbf{w}_{j}^k \gets \mathbf{w}_{j}^k - \frac{\rho}{B} \sum_{m \in \mathcal{B}}\left( \bm{\alpha} + (\bm \lambda^\star)^{(m)} \right)^\intercal \nabla_{\mathbf{w}_j^{k}} \mathbf{u}^{(m)}, \ \forall j \in \mathcal{J}^k, \forall k \in \mathcal{K}$\;
    }
    }
    Output ${\mathbf h}$.
  \caption{NN training}
\end{algorithm}

\cg{As described in Algorithm \ref{algo_1}, the gradient for NN training is computed using the optimal dual variable \(\bm{\lambda}^\star\). In other words, the training process involves solving a dual optimization problem, which may reduce training efficiency. However, the training is conducted offline and performed only once. During online decision-making, the trained NN can efficiently predict cost-effective dispatch decisions within a very short time. An alternative to solving optimization problems during training is dual learning  \cite{Park_Van_Hentenryck_2023}, where a separate NN is trained to directly predict optimal dual variables. While this approach eliminates the need for solving optimization problems during training, it may introduce additional prediction errors that potentially degrade the quality of dispatch decisions. Given the offline nature of our training process and the importance of ensuring high-quality dispatch decisions, we choose to retain the the dual problem solving during NN training.}

\subsection{Generalization performance guarantee}
This section provides a theoretical generalization performance guarantee, which characterizes the suboptimality of the proposed neural RLD with a large probability. Specifically, we first define the \emph{generalization cost} of a hypothesis $\mathbf{h}$, denoted by $L(\mathbf h)$, as follows:
\begin{align}
L(\mathbf h) = \mathbb{E}_{(\mathbf x, \mathbf d) \sim \mathcal{Z}} \left[\ell(\mathbf h(\mathbf x), \mathbf d)\right]. \label{eqn_generalization_cost}
\end{align}
This generalization cost is the same as the objective in \eqref{eqn_optimal}, which represents the expected cost of the decision given by $\mathbf{h}$ over the true (and unknown) population distribution $\mathcal{Z}$. For any hypothesis $\mathbf h$, we are also interested how its generalization cost compared to that of the optimal one $\mathbf h^\star$:
\begin{align}
    \Delta L (\mathbf{h}) = {L}(\mathbf{h}) - {L}(\mathbf h^\star), \label{eqn_GEC}
\end{align}
where $ \Delta L (\mathbf{h})$ is referred to as the \emph{excess cost of $\mathbf h$}. 

As introduced in Section \ref{sec_E2E}, we train a L2-regularized NN $\widehat{\mathbf h}$ to act as our hypothesis. In order to calculate the excess cost of $\widehat{\mathbf h}$, two mild assumptions are in force:
\begin{itemize}
    \item \textbf{A1}: Data points in dataset $\mathcal{S}$ are identically and independently distributed (i.i.d.), and follow an unknown distribution $\mathcal{Z}$.
    \item \textbf{A2}: The expectation of day-ahead auxiliary information vector $\mathbf x$ is bounded, i.e., $\mathbb{E}_{\mathbf x} \left(\Vert \mathbf x \Vert_2^2\right) \leq (X^\mathrm{max})^2$ for some finite $X^{\max}$.
\end{itemize}
While Assumption {\bf A1} is commonly used in the prediction step of the predict-then-optimize paradigm (e.g., in probabilistic forecasting settings), especially when theoretical performance guarantees are developed, the precise form of it may not be satisfied in many application settings. This does not fundamentally limit the applicability of our results for two reasons: (a) In any applications where past data is used to guide future decisions, some forms of stationarity assumption is necessary (so past data indeed contain relevant information for the decision). Our results can be extended to settings with weaker forms of stationarity assumption leveraging statistical learning theory tools for non-i.i.d. data \cite{pmlr-v63-Gao09}. (b) When the historical dataset is not stationary as it is, there may be transformations (e.g., de-trending and differencing as done in time series analysis) to remove non-stationarity.

We next derive a \emph{probably approximately correct} (PAC) upper bound for the excess cost of $\widehat{\mathbf h}$, i.e., a bound for the suboptimality of $\widehat{\mathbf h}$ that holds with high probability. This is achieved by applying statistical learning theory tools (i.e., uniform convergence) to our neural RLD hypothesis class and loss function, for which we develop our own Lipschitz continuity bound and Rademacher complexity bound,\footnote{Rademacher complexity measures the richness of a hypothesis class. A hypothesis class with a high Rademacher complexity in general is more expressive but requires more data points to learn (e.g., with more parameters).} respectively. The final result is stated as follows:
\begin{theorem}[PAC bound of excess cost] \label{prop_2}
When Assumptions \textbf{A1} and \textbf{A2} hold, with probability at least $1-\delta$ for any small $\delta \in (0,1)$, the excess cost of hypothesis $\widehat{\mathbf h}$ is bounded by:
\begin{align}
\Delta L(\widehat{\mathbf h}) \leq \frac{4(2{W^{\max}})^{K-\frac{1}{2}} C_{\ell} X^\mathrm{max} + \sqrt{2\ln(2/\delta)}}{\sqrt{M}},
 \label{eqn_EC_bound_final}
\end{align}
where $C_{\ell}$ is a finite constant related to the Lipschitz continuity of the RLD loss $\ell$.
\end{theorem}

\emph{Proof}: Appendix \ref{app_1} provides the detailed steps for deriving \eqref{eqn_EC_bound_final} and the expression for $C_{\ell}$.

By \eqref{eqn_GEC}, the PAC bound provided by Theorem \ref{prop_2} can be regarded as the expected optimality gap of the hypothesis $\widehat{\mathbf h}$. Thus, we also refer to Theorem \ref{prop_2} as the \emph{generalization performance guarantee} of hypothesis $\widehat{\mathbf h}$. This PAC bound decreases to zero as more samples are used, which indicates that the performance of $\widehat{\mathbf h}$ gradually approaches to the optimal hypothesis $\mathbf h^\star$ with probability $1-\delta$ even for arbitrarily small $\delta >0$. Hence, we have the following corollary.
\begin{colloary}[Asymptotic optimality] \label{remark_3}
The neural RLD hypothesis $\widehat{\mathbf h}$ is asymptotically optimal, i.e., with probability at least $1-\delta$, $\Delta L(\widehat{\mathbf h}) \to 0 $ as $M \to \infty$ for any $\delta \in (0,1)$.
\end{colloary}

\cg{Although the neural RLD method and the generalization performance guarantee are established for L2-regularized NNs, they are fully compatible with other types of NNs, including L1-regularized NNs, recurrent neural networks (RNNs), and convolutional neural networks (CNNs). First, regardless of the type of NN employed, the model can always be trained to predict cost-effective first-stage decisions \(\mathbf{u}\) based on observable features \(\mathbf{x}\). Second, extending the performance guarantee to other types of NNs requires only minor modifications to the derivation process. Detailed descriptions of these modifications are provided in Appendix \ref{app_1}.}

Note the RLD formulation in this paper is based on the DC power flow model without considering some realistic settings, such as AC power flow formulation, ramping constraints, and production limits, while allowing unrestricted generation reduction. Nevertheless, the proposed neural RLD method can be extended to incorporate these settings. A detailed discussion on these extensions is provided in Appendix \ref{app_4}.

\section{Case Study} \label{sec_case}
\subsection{Simulation setup}
\subsubsection{Test systems}
\begin{figure}[h]
	\vspace{-8mm}
	\subfigbottomskip=-4pt
	\subfigcapskip=-4pt
	\centering
	\subfigure[5-bus case]{\includegraphics[width=0.49\columnwidth]{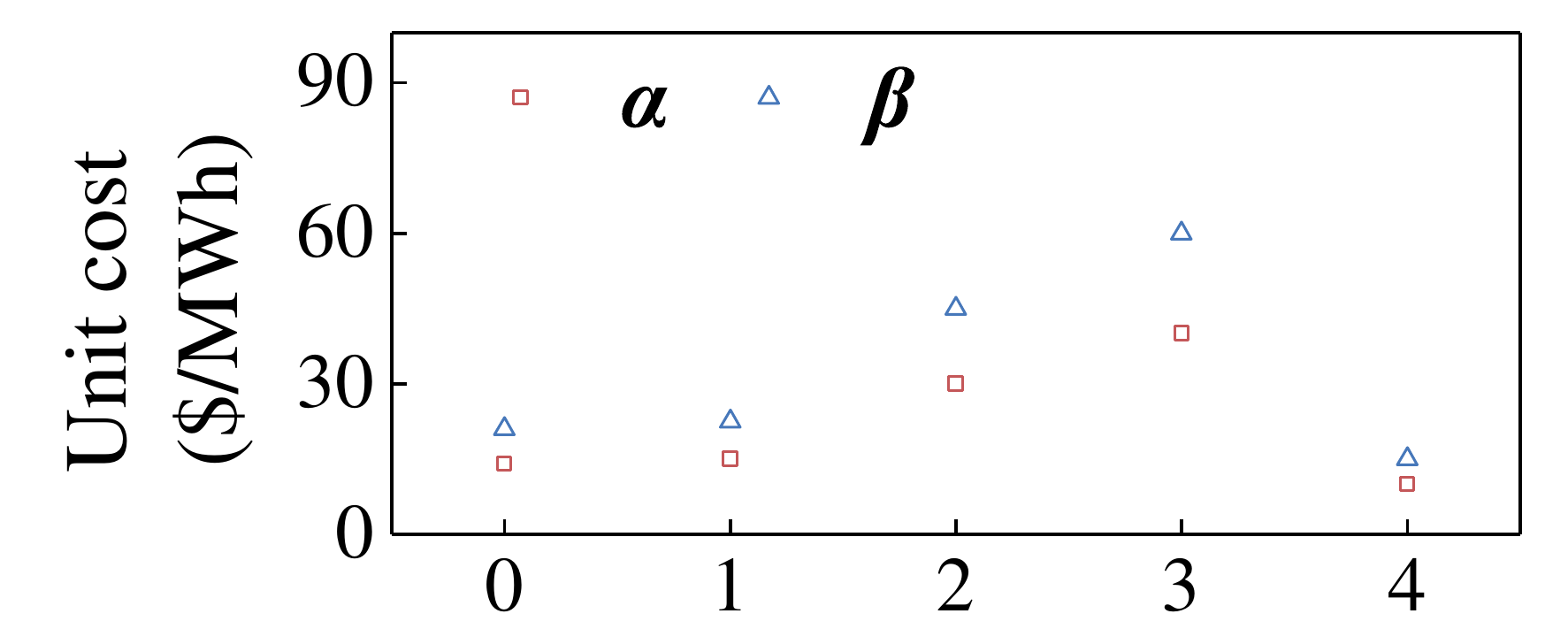}} 
	\subfigure[118-bus case]{\includegraphics[width=0.49\columnwidth]{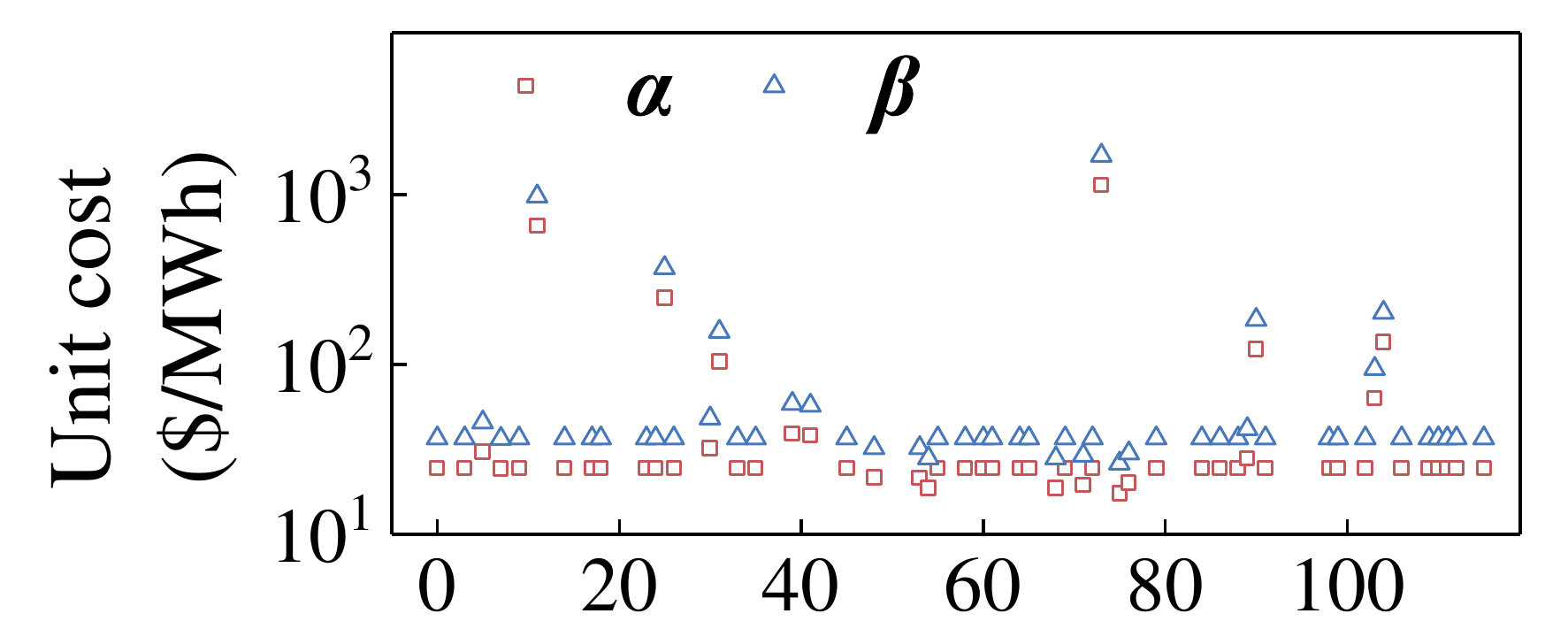}}
	\subfigure[300-bus case]{\includegraphics[width=0.49\columnwidth]{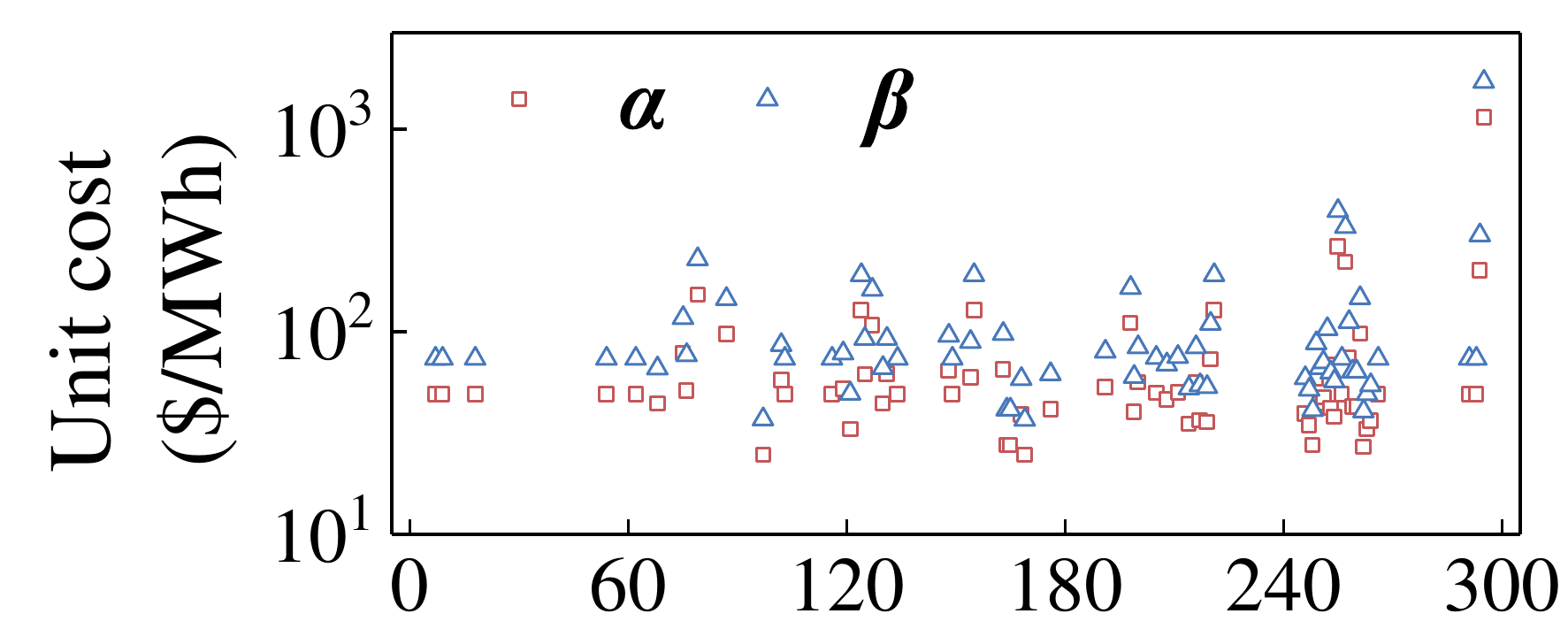}}
	\subfigure[1354-bus case]{\includegraphics[width=0.49\columnwidth]{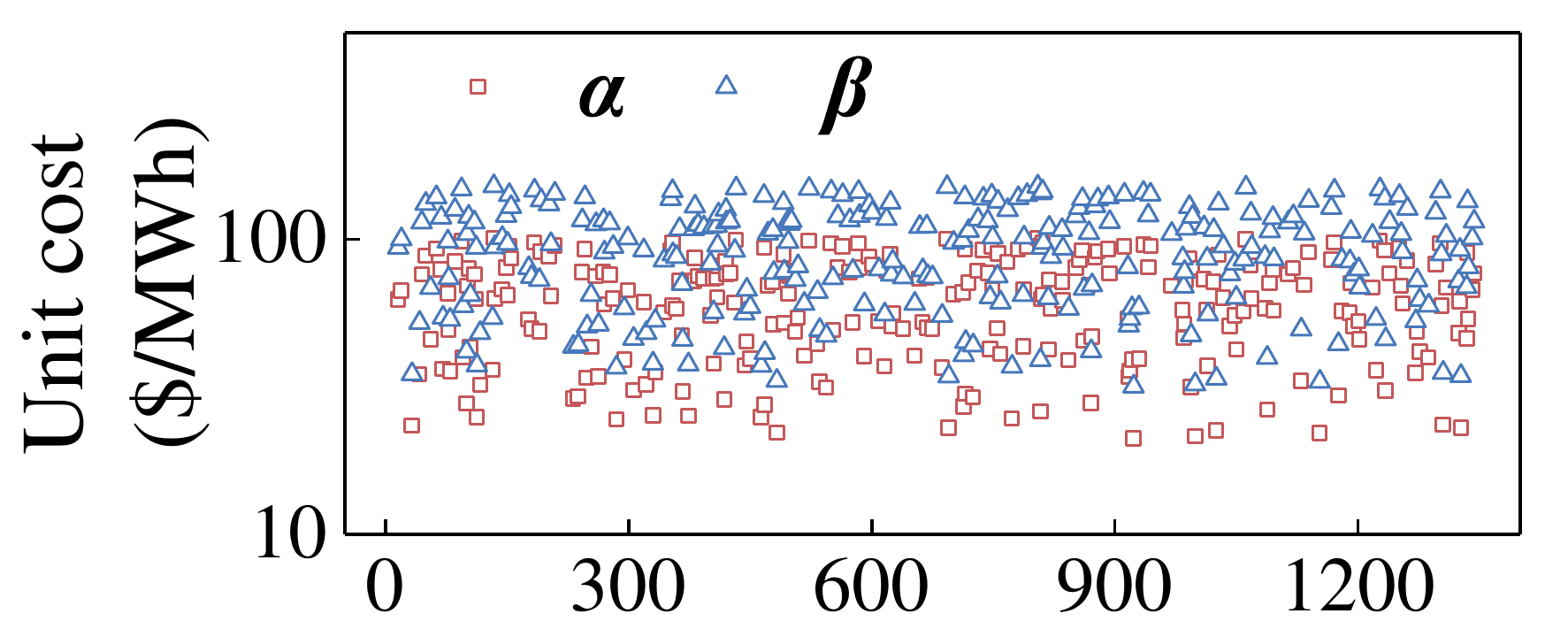}}
 	\caption{Unit generation costs used in different test cases. Note the first-stage unit generation cost $\bm \alpha$ is always smaller than the second-stage one $\bm \beta$.}
	\label{fig_unitcost}
	\vspace{-2mm}
\end{figure}

To validate the performance of the proposed framework, we implement our case studies based on four different systems, i.e, the IEEE 5-bus, 118-bus, 300-bus, and 1354-bus test systems. The numbers of buses, line, loads, and generators in different test systems are summarized in Table \ref{tab_parameter}. 
The unit costs in the first-stage and second-stage, i.e., $\bm \alpha$ and $\bm \beta$, are illustrated in Fig. \ref{fig_unitcost}. Other parameters remain the same as the standard IEEE test cases. 
\begin{table}[h]
\vspace{-2mm}
	\small
	\centering
\begin{threeparttable}
\caption{Descriptions of our test cases and structures of the corresponding NNs}\label{tab_parameter}
\begin{tabular}{ccccc}
\hline
\rule{0pt}{11pt}
  Cases & \makecell{Line \\ number} & \makecell{Load \\ number} & \makecell{Generator \\ number} & \makecell{Structures of \\ of NNs} \\ \hline
  5-bus &  6 & 5 & 5  & (5, 5, 5) \\
    118-bus & 173 & 99 & 54 & (20, 20, 20) \\
    300-bus &  411 &  199 &  69  & (40, 40, 40) \\
    1354-bus & 1991 & 673 & 260 & (60, 60, 60) \\
    \hline
\end{tabular}\end{threeparttable} 
\end{table}

\subsubsection{Benchmarks}
The following two benchmarks are introduced to verify the benefits of the proposed neural RLD:
\begin{itemize}
\item Two-step predict-then-optimize paradigm used in \cite{6496182} (termed ``\textsf{Two-step}"). It first follows \cite{8982039} and trains a mixture density network to learn the distribution of net demands. According to our test, this model shows high accuracy, e.g., the relative errors for predicting the mean and covariance are smaller than 3\%. Then, the scenario-based method is employed to solve problem \eqref{eqn_RLD_2nd}, in which 30 discrete scenarios are generated to approximate the predicted distribution. 
\item Imitation learning used in \cite{9585298} (termed ``\textsf{Imitation}''). This benchmark trains a NN using the MSE loss to directly learn the mapping from auxiliary information $\mathbf x$ to the dispatch decision $\mathbf u$. The training labels, i.e., historical optimal decisions, are generated by solving the deterministic dispatch problems given the realizations of $\mathbf d$ in the past delivery intervals.
\end{itemize}

Note that benchmark \textsf{Two-step} is a generalized version of the classic RLD method \cite{6818365}. The classic RLD assumes that future net demand follows a Gaussian distribution \cite{6818365}. In contrast, benchmark \textsf{Two-step} employs a mixture density network to predict the probability distribution of future net demand, which accommodates both Gaussian and non-Gaussian distributions and thereby enhancing prediction accuracy. As a result, the benchmark is expected to outperform the classic RLD.

\subsubsection{Data generation}
For each test case, we randomly generate 5,000 samples of $\mathbf x \in \mathbb R^{p}$ from a uniform distribution on the interval $[0, 1]^{p}$. Then, the uncertain net demand $\mathbf d \in \mathbb R^{\Nbus}$ is generated based on the following rule:\footnote{One may anticipate good prediction performance and perhaps superior performance of the \textsf{Two-step} paradigm due to the linearity between $\mathbf d$ and $\mathbf x$ in data generation. While the former observation is indeed the case, the latter turned out to be false due to the mismatch of the two steps in the predict-then-optimize paradigm. In fact, we choose the simple rule \eqref{eqn_data_rule} precisely to highlight the importance of end-to-end decision making. }
\begin{align}
\mathbf d = \bm \Omega \left(\mathbf{x} * (\bm 1 + 0.15 \cdot \bm \omega)\right), \label{eqn_data_rule}
\end{align}
where $\bm \Omega \in \mathbb R^{\Nbus \times p}$ is a fixed coefficient matrix; $*$ represents the element-wise multiplication. A random variable $\bm \omega \in \mathbb R^{p}$ is introduced to simulate the stochastic characteristic of net demands. Here, we assume that $\bm \omega$ follows the standard Gaussian distribution. {Among the 5,000 samples, 4,000 are used as the training set, while the remaining 1,000 are used as the testing set to evaluate the suboptimality of the different methods in our case study.}

\subsubsection{Structures of NNs}
The NNs used in all methods are composed of fully connected layers. Their structures are summarized in the column ``Structures of NNs" in Table \ref{tab_parameter}, where every number in brackets represents the neuron number in a hidden layer. For example, (5, 5, 5) represents the corresponding NN has three hidden layers, each with five neurons. {For larger power networks, we utilize more neurons in each layer of the NN as we expect a richer hypothesis class is needed. For example, consider the trivial case where $\mathbf x = \mathbf d$ and $\mathcal D(\mathbf x)$ is a point mass, i.e., the RLD problem reduces to a simple deterministic linear problem. In this case, $\mathbf{u}^\star(\mathbf{x})$ is a piecewise function determined by a parametric linear program \cite{still2018lectures}. The number of pieces corresponds to the number of \emph{critical regions}, each defined as a set of $\mathbf{x}$ over which the reduced linear problem exhibits the same constraint binding pattern at the solution $\mathbf{u}^\star(\mathbf{x})$. Thus, the number of critical regions depends on the number of possible constraint binding patterns. As larger networks have more inequality constraints (e.g., line flow constraints), the resulting mapping $\mathbf{u}^\star(\mathbf{x})$ is necessarily more complex and requires more neurons to approximate. In standard machine learning terminology, the process of determining the number of neurons is referred to as \emph{model selection} and can be done via, e.g., cross-validation. }


\subsection{Performance evaluation}
\subsubsection{Convergence} 

\begin{figure}[h]
	\vspace{-4mm}
	\subfigbottomskip=-4pt
	\subfigcapskip=-4pt
	\centering
	\subfigure[5-bus case]{\includegraphics[width=0.49\columnwidth]{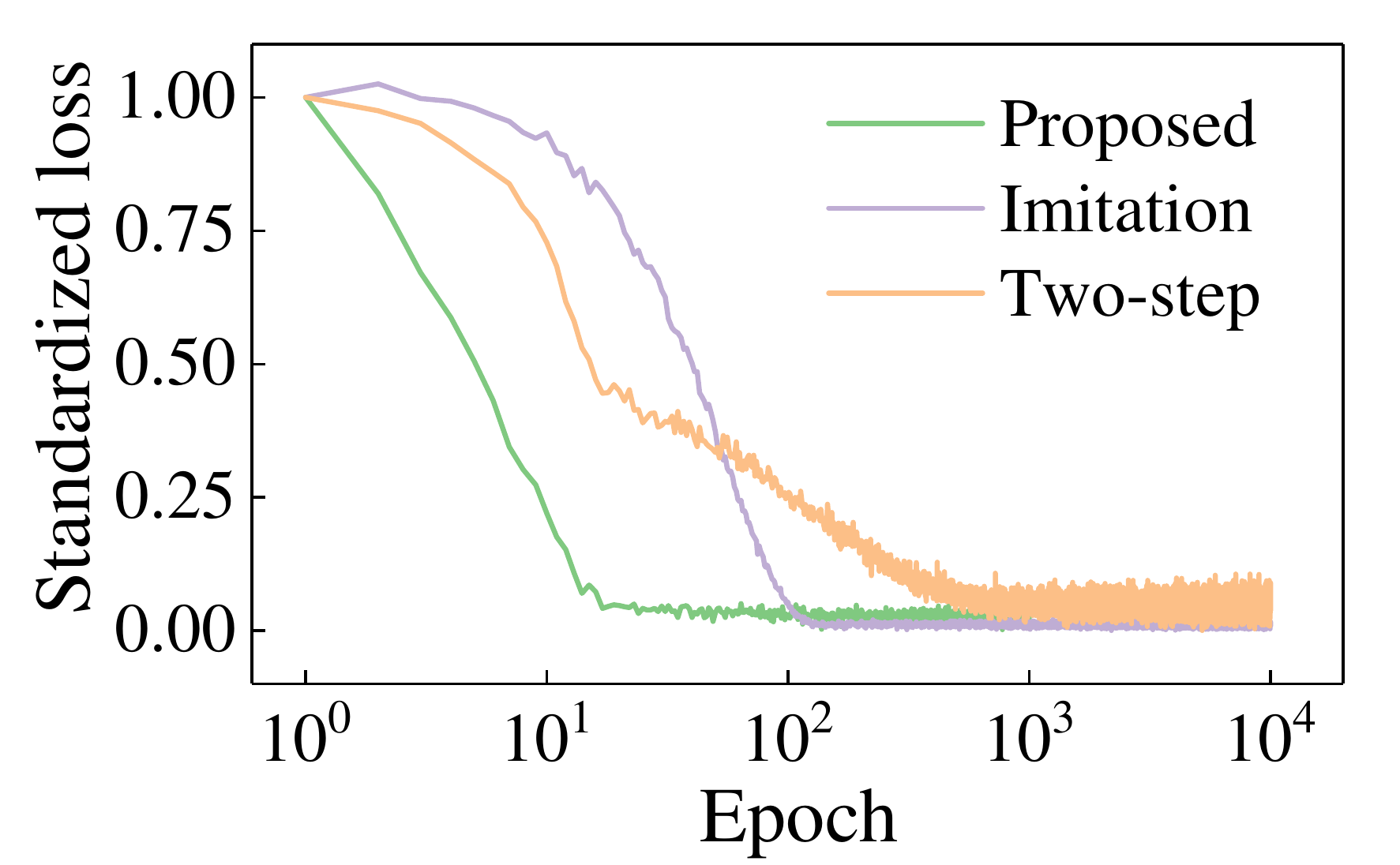}} 
	\subfigure[118-bus case]{\includegraphics[width=0.49\columnwidth]{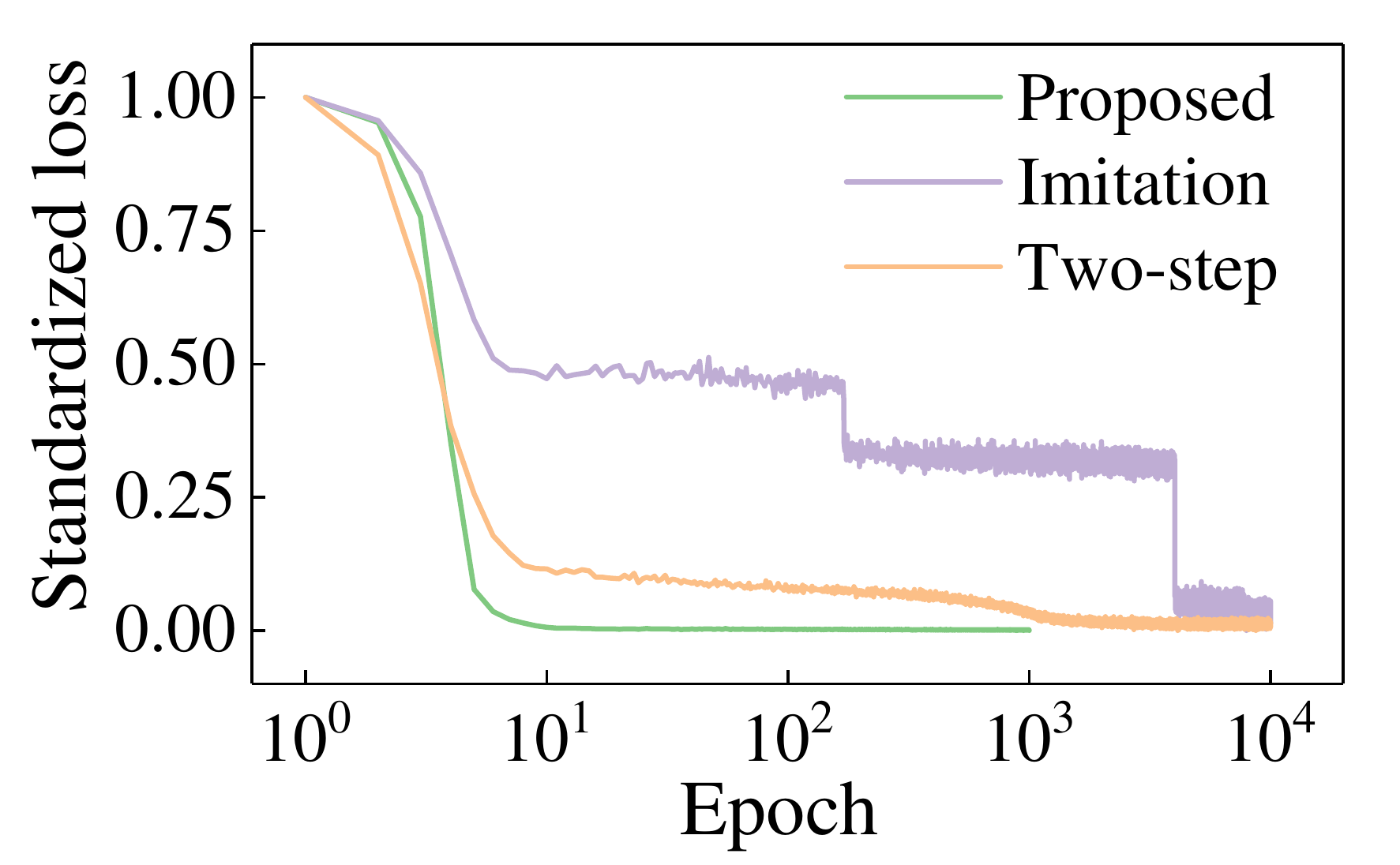}}
	\subfigure[300-bus case]{\includegraphics[width=0.49\columnwidth]{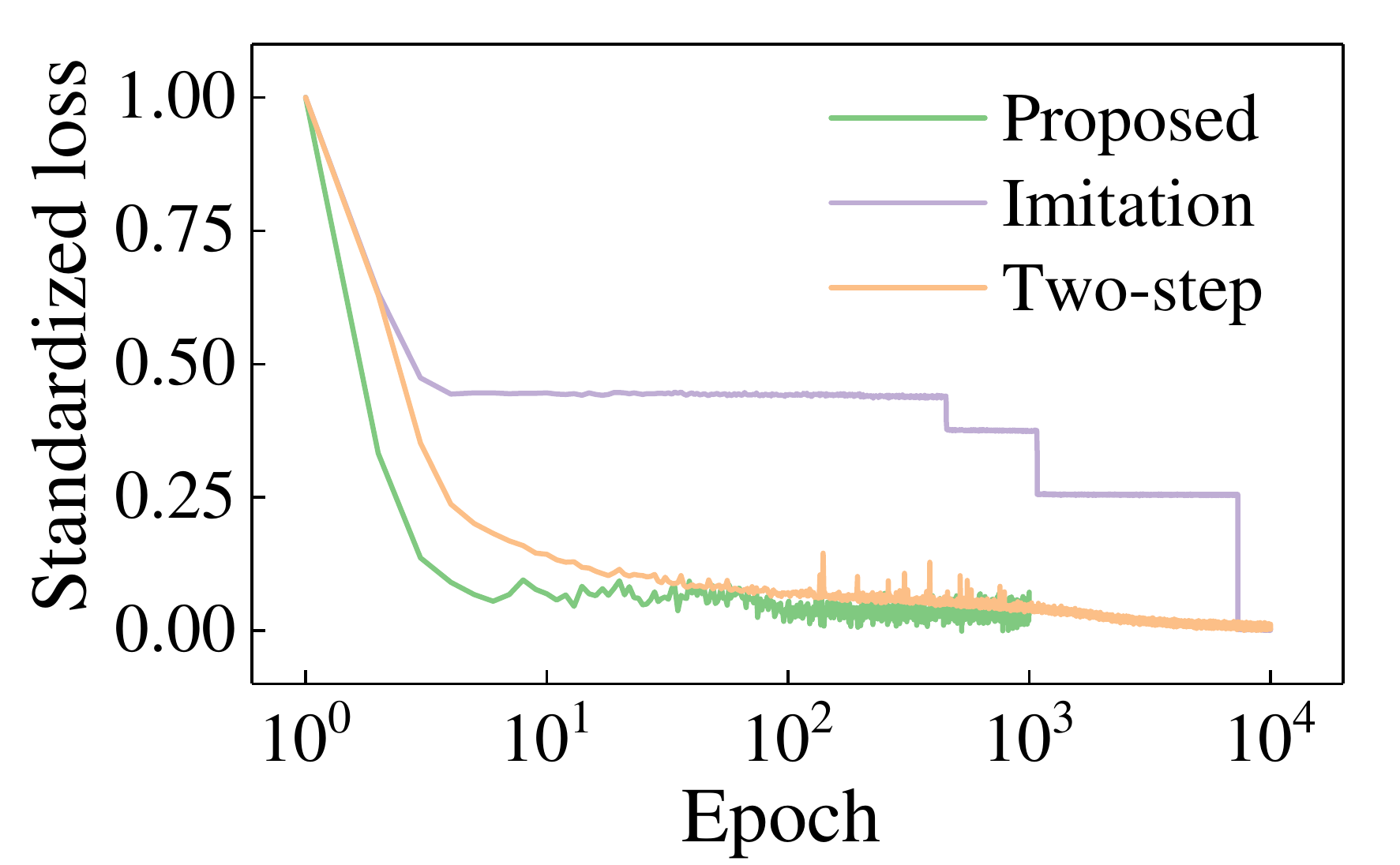}}
	\subfigure[1354-bus case]{\includegraphics[width=0.49\columnwidth]{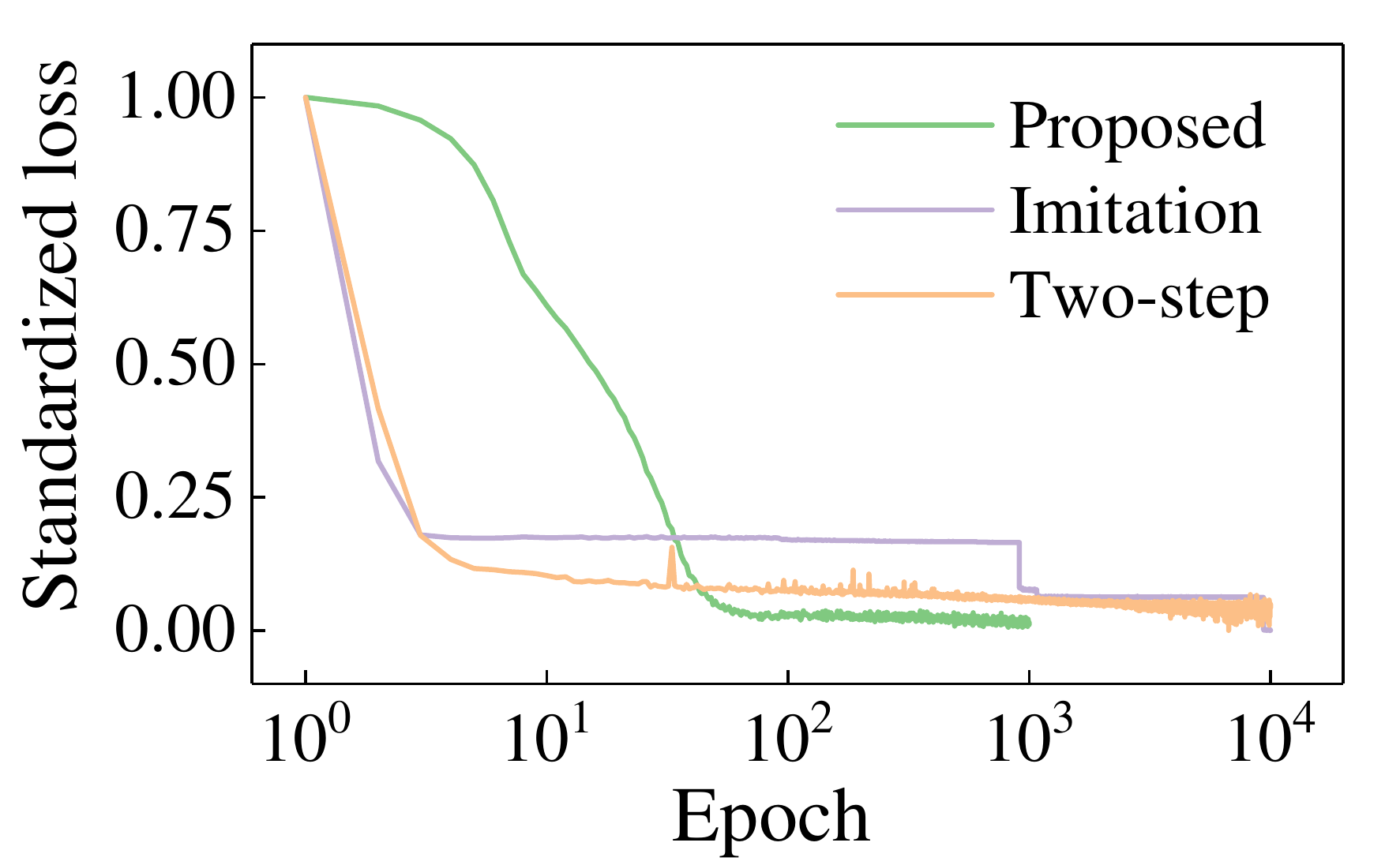}}
 	\caption{Standardized training losses of the NNs trained by different methods.}
	\label{fig_trainingloss}
\end{figure}

Fig. \ref{fig_trainingloss} demonstrates the training losses of the NNs in each test case. Note all the training losses are standardized to [0,1] for a fair comparison. In every case, the proposed method can converge in only 100 epochs. In contrast, both the benchmarks may not converge even after 1,000 epochs.\footnote{
The convergence of our proposed method reaches a saturation point (i.e., starts to slow down) with fewer epochs in the first three cases. However, in the 1354-bus case, it appears later than the other two methods. Nevertheless, an earlier saturation point does not inherently indicate better training efficiency or lower suboptimality. For instance, although the saturation point of the proposed method occurs later in the 1354-bus case, its standardized loss nearly stabilizes after 100 epochs. In contrast, the loss curve of benchmark \textsf{Imitation} shows a significant drop around 1000 epochs, while that of benchmark \textsf{Two-step} continues to decrease even beyond 1000 epochs. These variations in the timing of the saturation point across different methods and test cases may result from the complex interplay among various factors, such as the NN's scale and the definition of the loss function.
} 
According to Proposition \ref{prop_1}, the gradient of the RLD loss is well-defined, so the proposed model can be efficiently trained by the SGD. In benchmark \textsf{Imitation}, the NN is trained using the MSE to learn the mapping from auxiliary information $\mathbf{x}$ to the optimal decision $\mathbf{u}^\star$. Since the RLD problem \eqref{eqn_RLD} is not strictly convex, this mapping may be discontinuous \cite{elmachtoub2022smart}, which is hard to learn. Benchmark \textsf{Two-step} trains a mixture density network to learn the distribution of net demands for a given $\mathbf{x}$ \cite{8982039}. This mixture density network needs to predict the means and covariances for multiple Gaussian components, while a set of parameters need to be learned for each Gaussian component. Thus, it has more learnable parameters than the common neuron network even if they have the same hidden layers and neurons, leading to lower training efficiency.

\subsubsection{Suboptimality}
\begin{figure}[h]
	\vspace{-6mm}
	\subfigbottomskip=-4pt
	\subfigcapskip=-4pt
	\centering
	\subfigure[5-bus case]{\includegraphics[width=0.49\columnwidth]{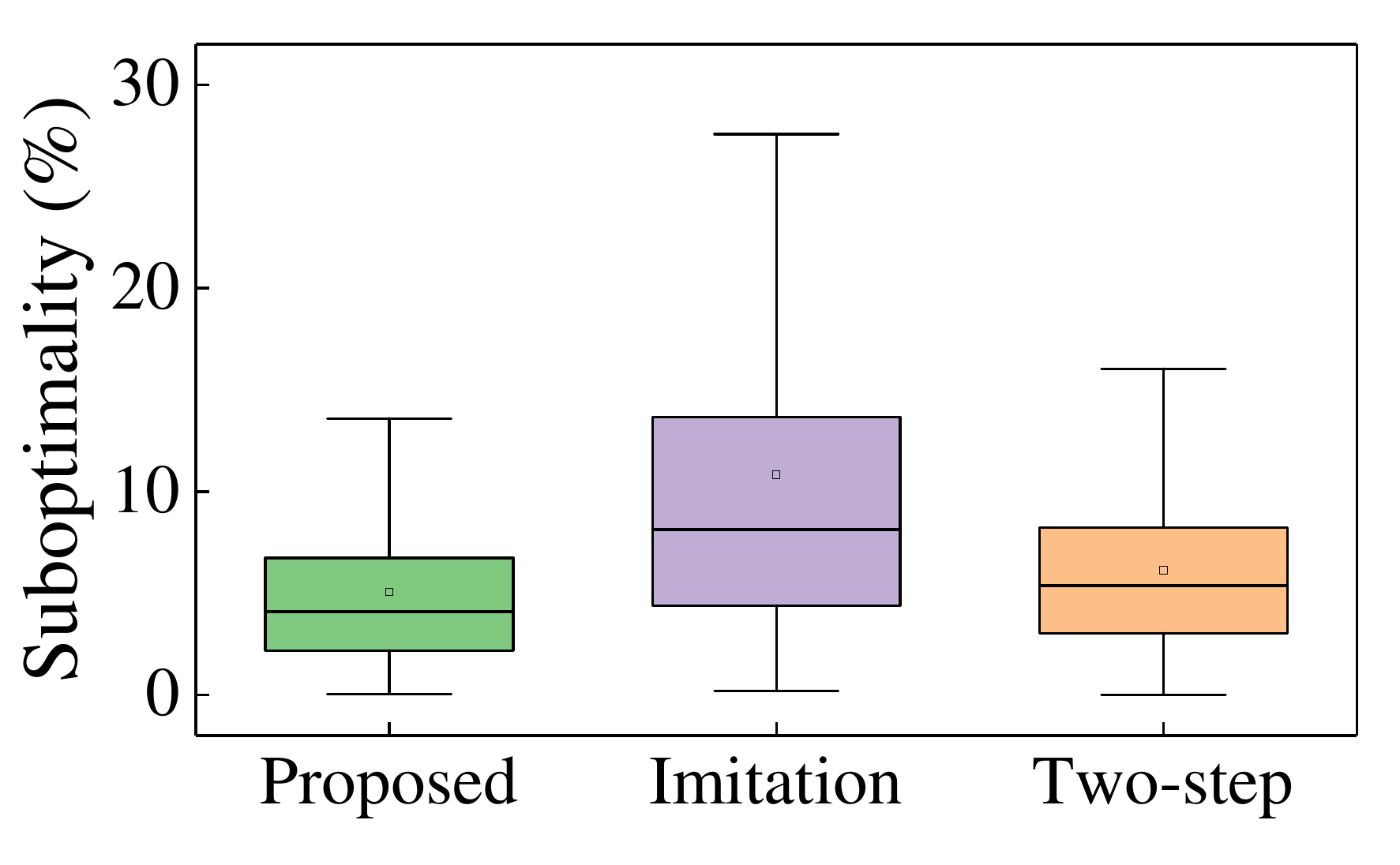}} 
	\subfigure[118-bus case]{\includegraphics[width=0.49\columnwidth]{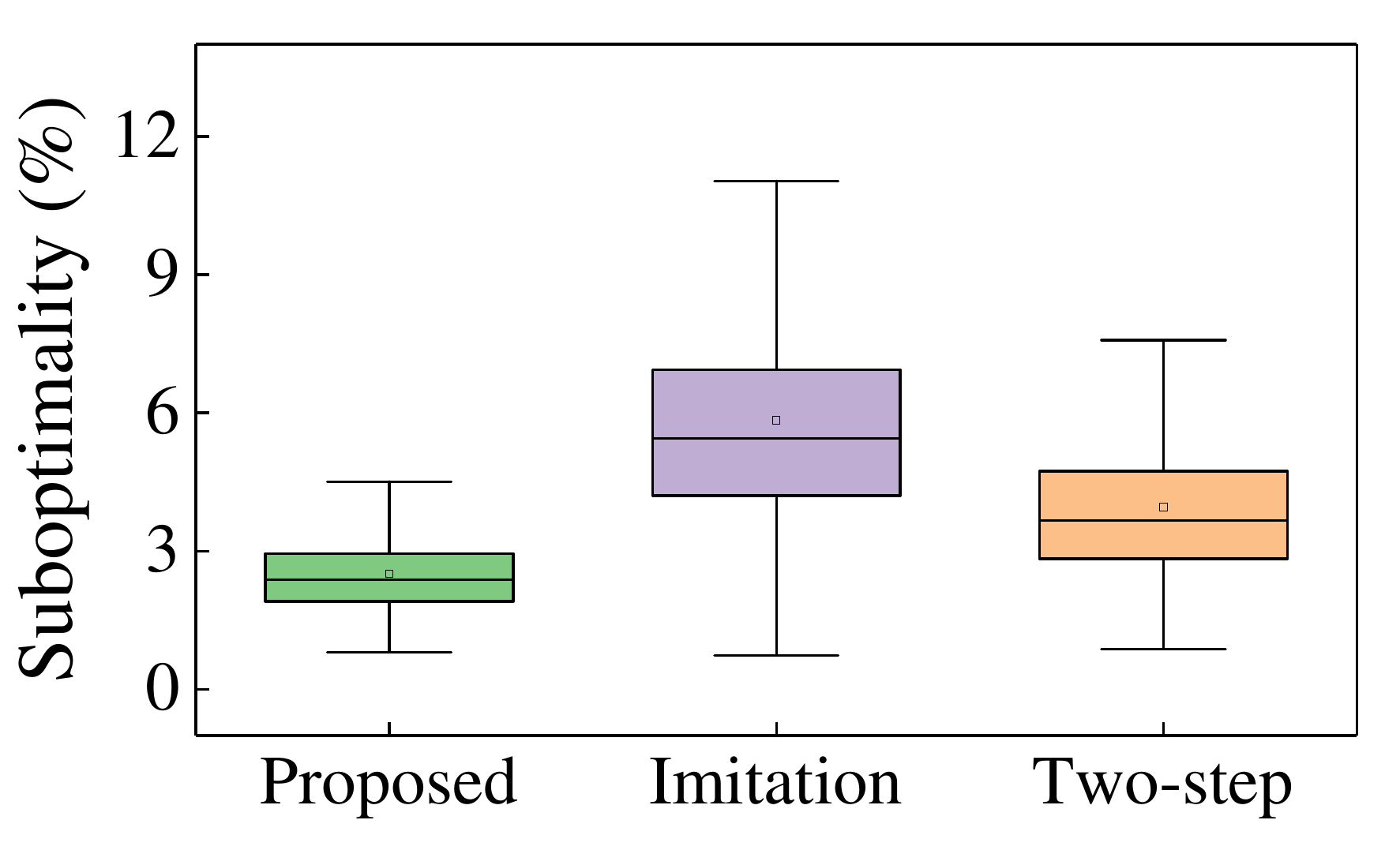}}
	\subfigure[300-bus case]{\includegraphics[width=0.49\columnwidth]{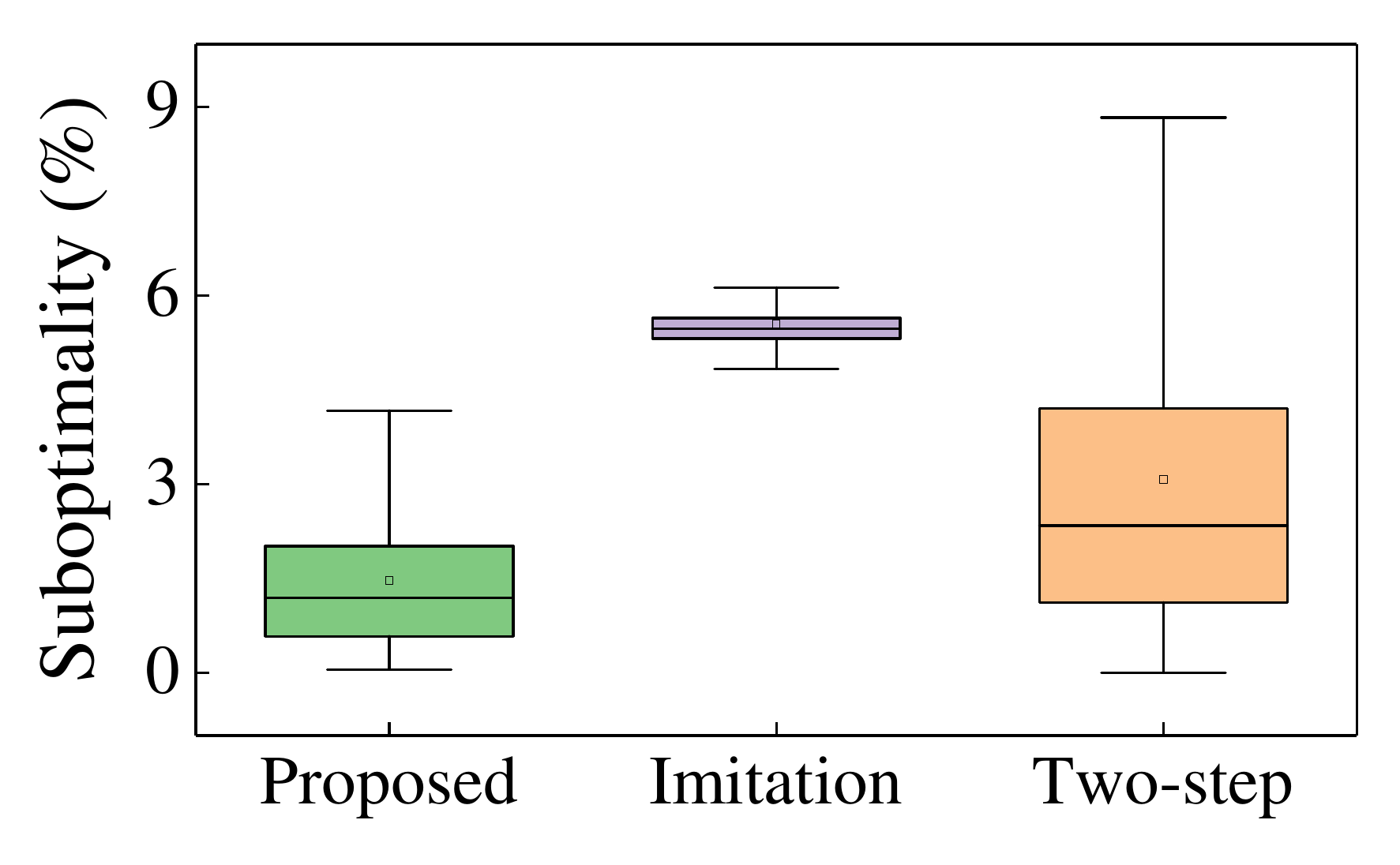}}
	\subfigure[1354-bus case]{\includegraphics[width=0.49\columnwidth]{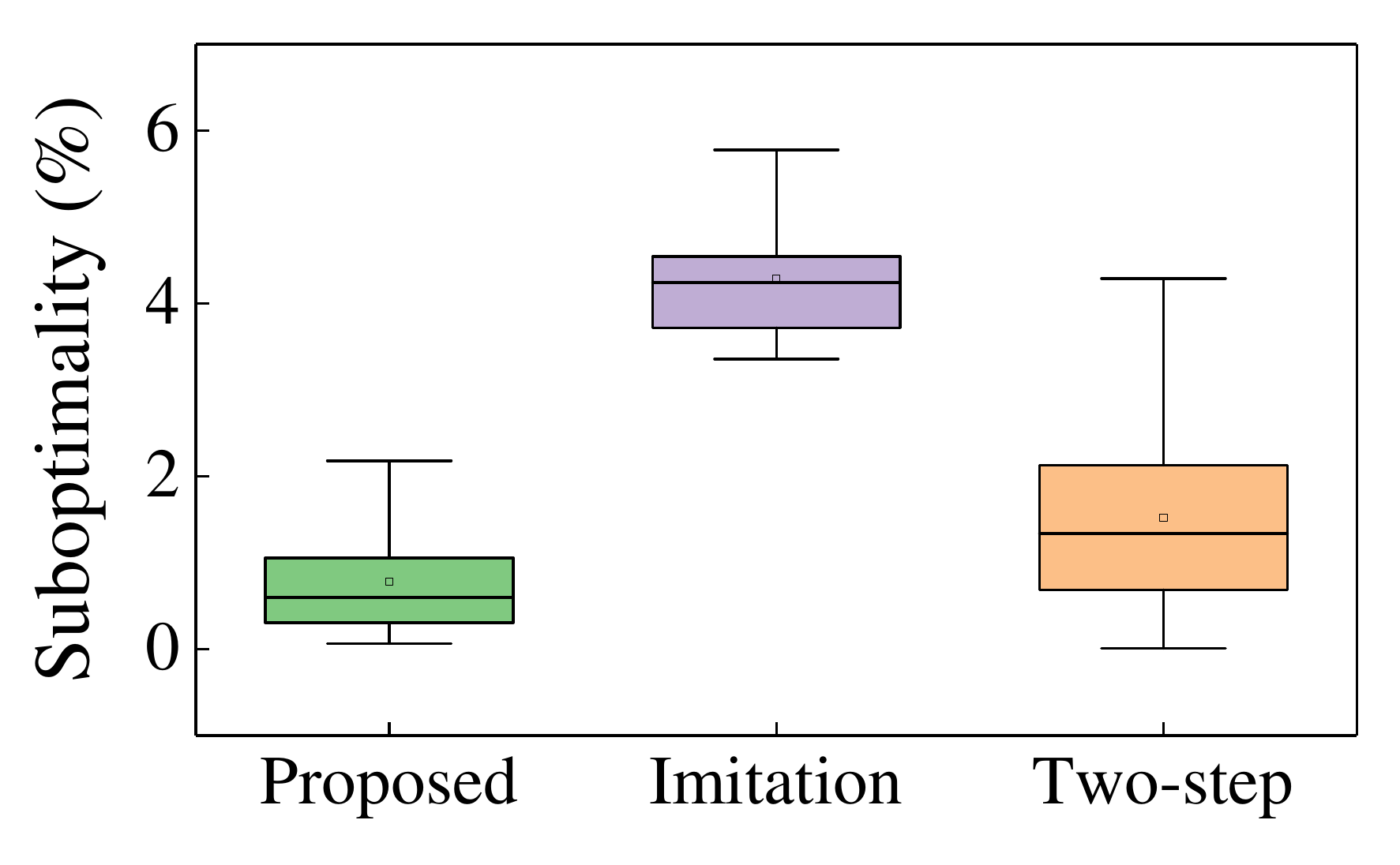}}
 	\caption{Normalized suboptimality given by different methods.}
	\label{fig_regret}
	\vspace{-2mm}
\end{figure}

We compare the normalized suboptimality of different methods, which is defined as:
\begin{align}
\mathrm{Suboptimality} = \frac{\ell(\widehat{\mathbf u}, \mathbf d) - \ell({\mathbf u^\star}, \mathbf d)}{\ell({\mathbf u^\star}, \mathbf d)}, 
\end{align}
where $\widehat{\mathbf u}$ represents the dispatch decision given by each method. The numerator represents the excess generation cost between $\widehat{\mathbf u}$ and the optimal decision ${\mathbf u^\star}$, while the denominator is the optimal objective of \eqref{eqn_RLD_2nd}. We calculate this normalized suboptimality for every data point in the testing set. The results are given in Fig. \ref{fig_regret}. In all test cases, the normalized suboptimality of the proposed method is significant smaller than that of the two benchmarks. For instance, the average suboptimality of the proposed method is only 1.46\% in the 300-bus test case, while this value is 5.54\% and 3.07\% in benchmarks \textsf{Imitation} and \textsf{Two-step}, respectively. As aforementioned, the RLD loss can more accurately reflect the suboptimality of decisions compared to the MSE loss used in benchmark \textsf{Imitation}. {Moreover, the proposed method is an end-to-end approach that merges the prediction and optimization steps. It can bypass the ``mismatch" issue encountered in benchmark \textsf{Two-step}, where prediction and optimization are considered separately. Therefore, it demonstrates the lowest suboptimality.}

\subsubsection{Computational efficiency}
\begin{figure}[h]
	\vspace{-4mm}
	\subfigbottomskip=-4pt
	\subfigcapskip=-4pt
	\centering
	\subfigure[5-bus case]{\includegraphics[width=0.49\columnwidth]{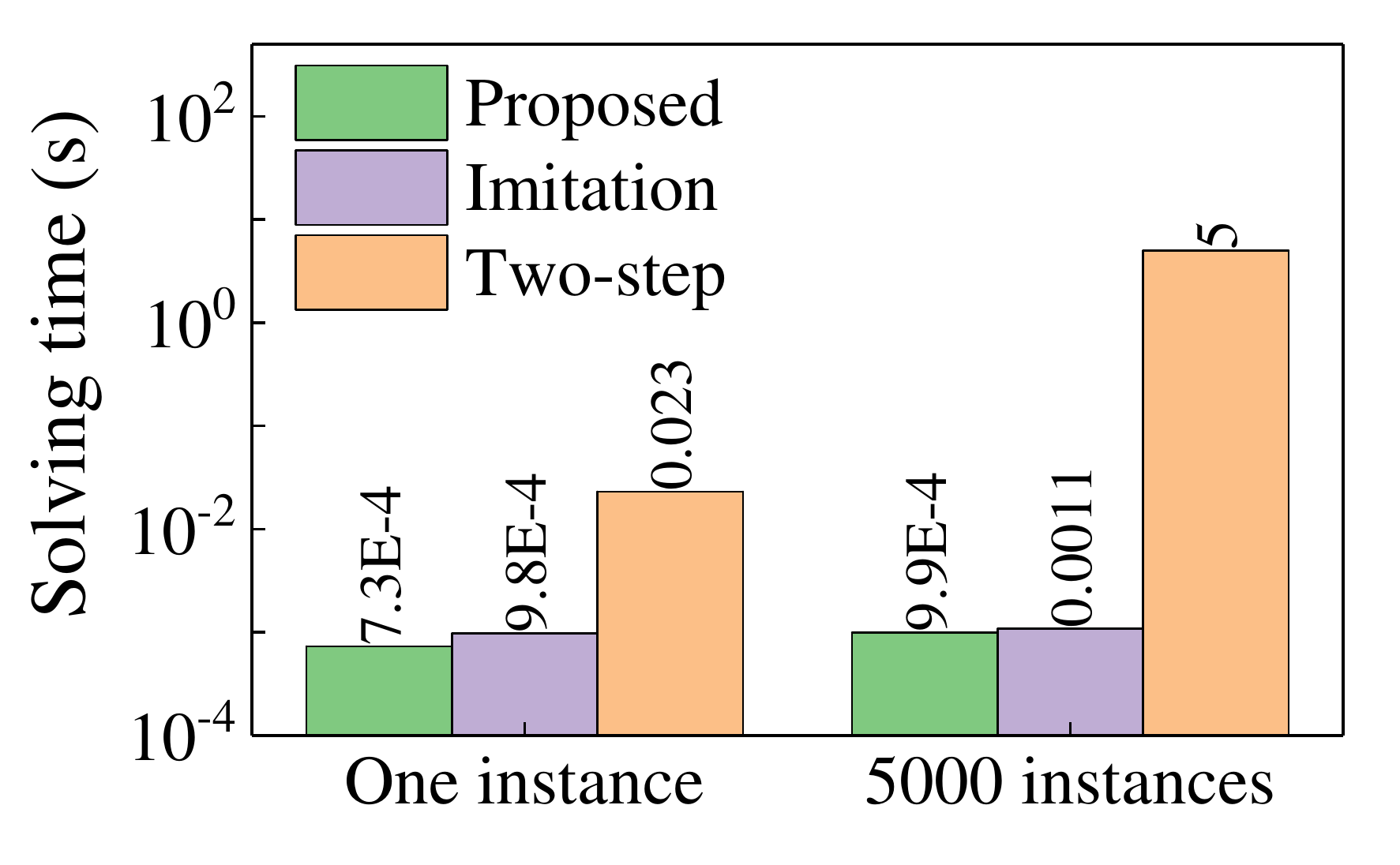}} 
	\subfigure[118-bus case]{\includegraphics[width=0.49\columnwidth]{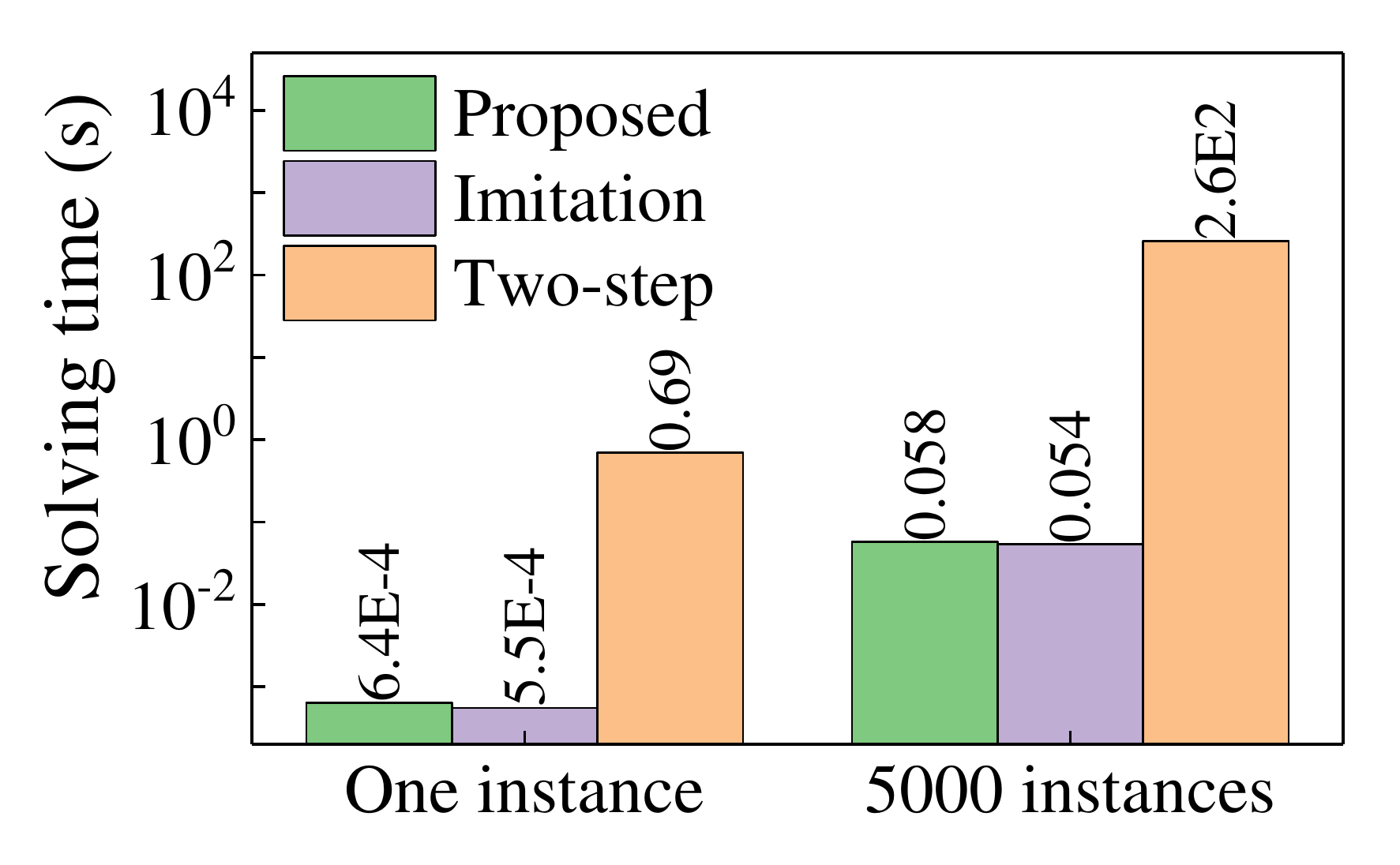}}
	\subfigure[300-bus case]{\includegraphics[width=0.49\columnwidth]{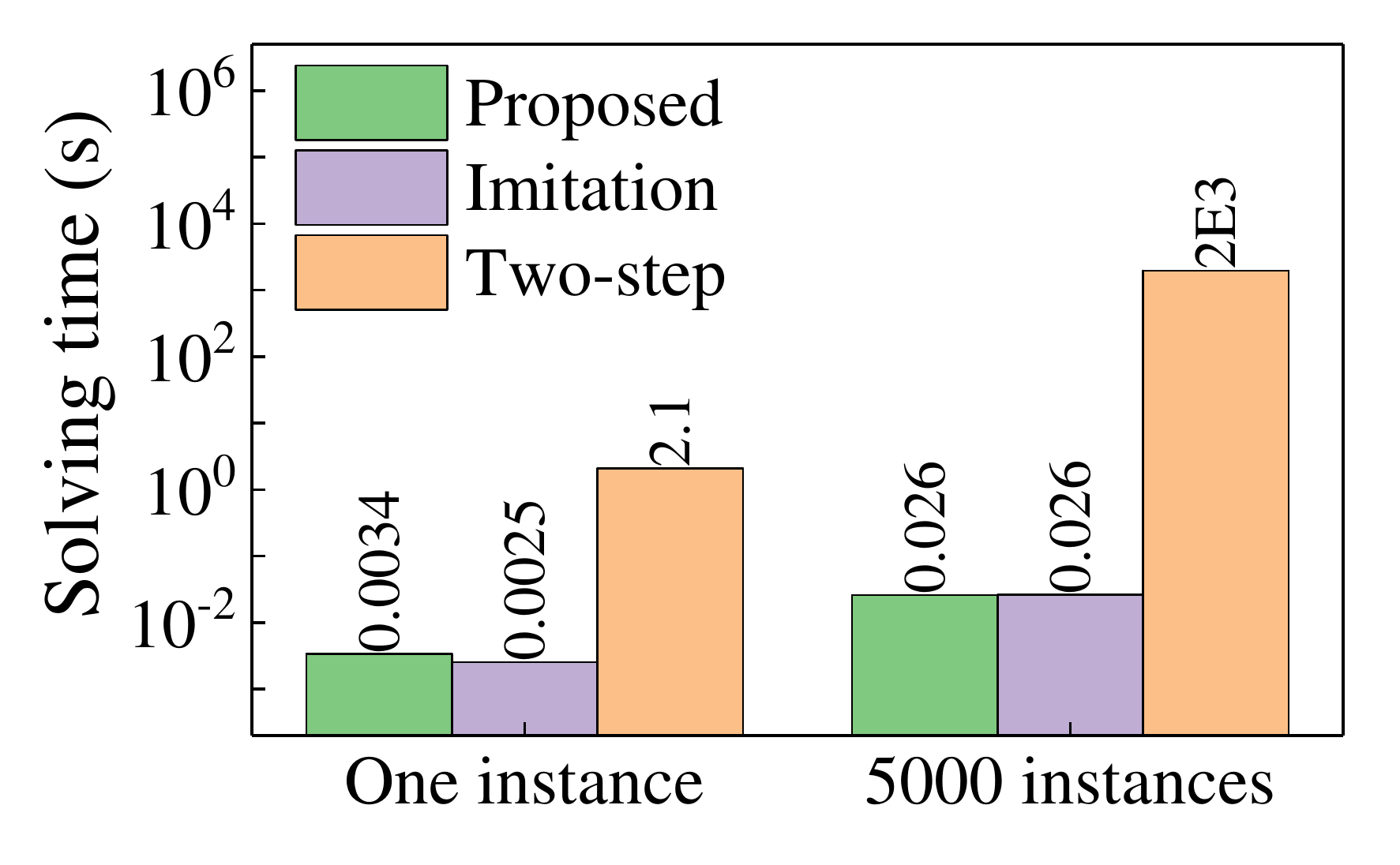}}
	\subfigure[1354-bus case]{\includegraphics[width=0.49\columnwidth]{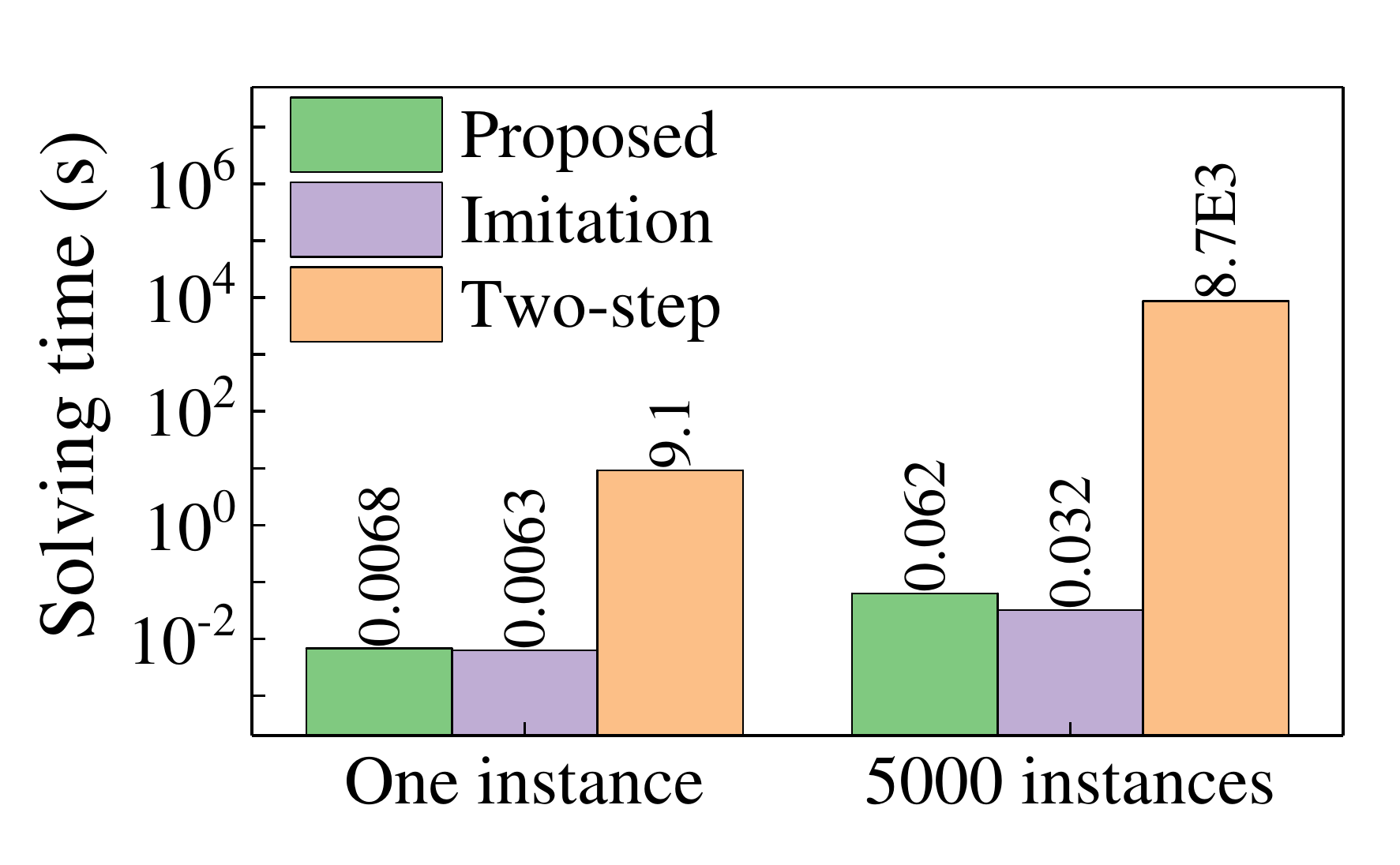}}
 	\caption{Times for solving one and 5,000 instances with different methods.}
	\label{fig_solvingTime}
\vspace{-4mm}
\end{figure}

Fig. \ref{fig_solvingTime} compares the required time for solving one single instance and 5000 instances of \eqref{eqn_RLD} using different methods, where benchmark \textsf{Two-step} is solved in parallel with 30 threads. {Since both the proposed neural RLD and benchmark \textsf{Imitation} replace the solving process with the forward pass of NNs, they can generate dispatch decisions in a very short time. In contrast, the \textsf{Two-step} benchmark still requires solving a stochastic program after predicting the distribution of net demand, resulting in significantly longer computation time. For instance, the solving times for the proposed method and the \textsf{Imitation} benchmark are at least two orders of magnitude smaller than those of benchmark \textsf{Two-step}.} 

\subsection{Sensitivity analysis}
\subsubsection{Neuron number}
\begin{figure}[h]
		\vspace{-2mm}
	\centering
	{\includegraphics[width=1\columnwidth]{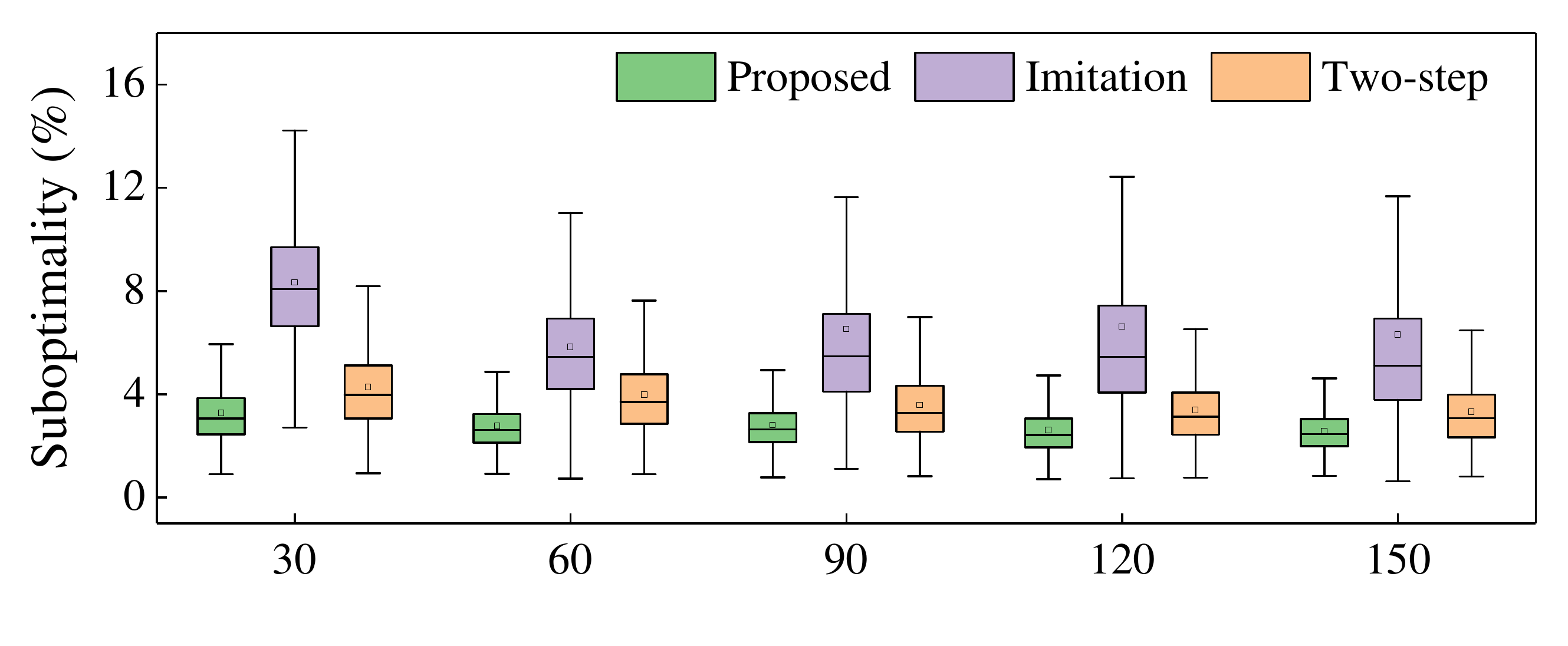}}
	\vspace{-10mm}
 	\caption{Normalized suboptimality of all methods under different neuron numbers in the 118-bus test case.}
	\label{fig_regret_neuronNum}
	\vspace{-2mm}
\end{figure}

\begin{figure}[h]
		\vspace{-4mm}
	\centering
	{\includegraphics[width=1\columnwidth]{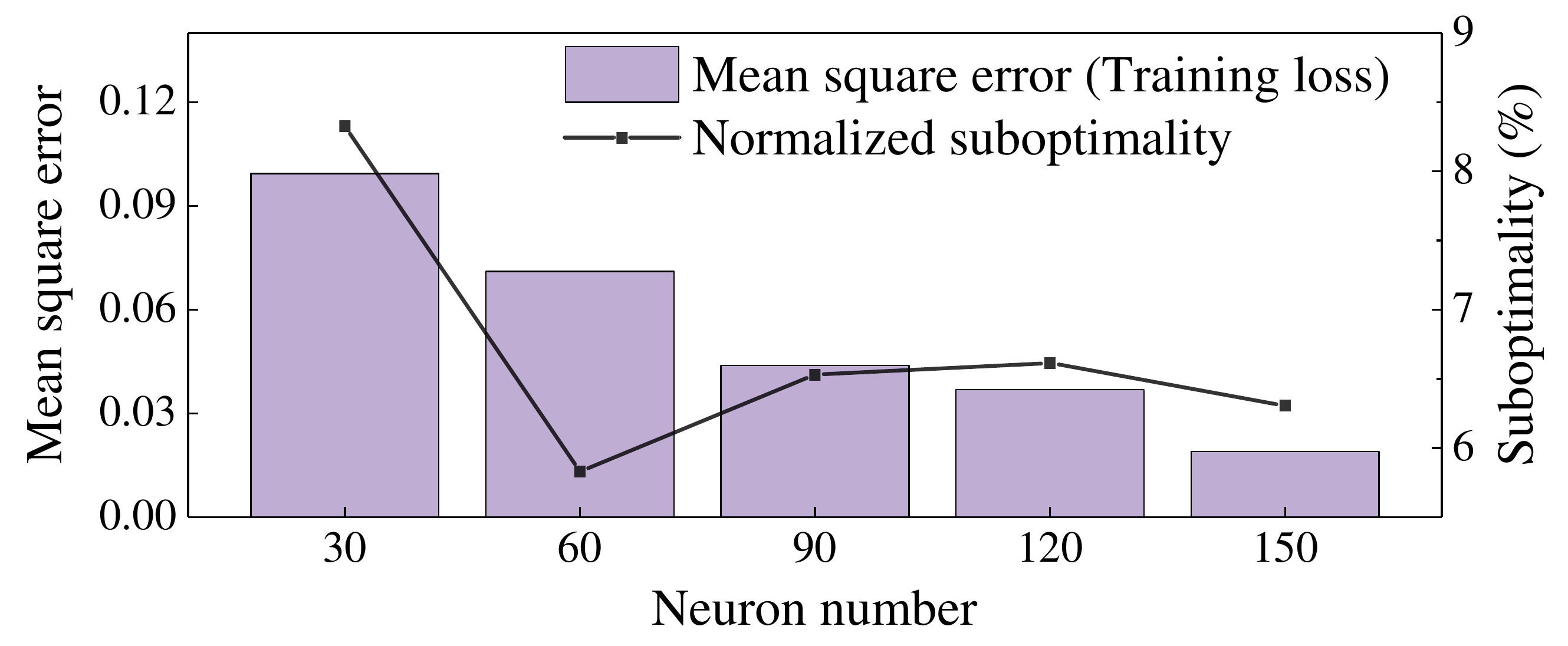}}
	\vspace{-8mm}
 	\caption{Training losses, i.e., the MSE, and average normalized suboptimality of the imitation learning under different neuron numbers. Obviously, a smaller MSE may not always result in lower suboptimality. }
	\label{fig_MSE_neuronNum}
	\vspace{-4mm}
\end{figure}

We investigate the impacts of neuron numbers on the normalized suboptimality of different methods based on the 118-bus test system. The results are illustrated in Fig. \ref{fig_regret_neuronNum}. The proposed method consistently exhibits the lowest suboptimality in all cases. Meanwhile, even with only a few neurons, the proposed method demonstrates low suboptimality. This superiority results from the RLD loss, enabling more accurate evaluation on the suboptimality of decisions compared to the MSE. This conclusion can be further validated by Fig. \ref{fig_MSE_neuronNum}, in which the MSE values and average suboptimality of benchmark \textsf{Imitation} method are illustrated. While increasing neuron number can reduce the MSE, the suboptimality does not exhibit a consistent decrease. 

\subsubsection{Hidden layer number}

\begin{figure}[h]
	\vspace{-4mm}
	\centering

 	\includegraphics[width=1\columnwidth]{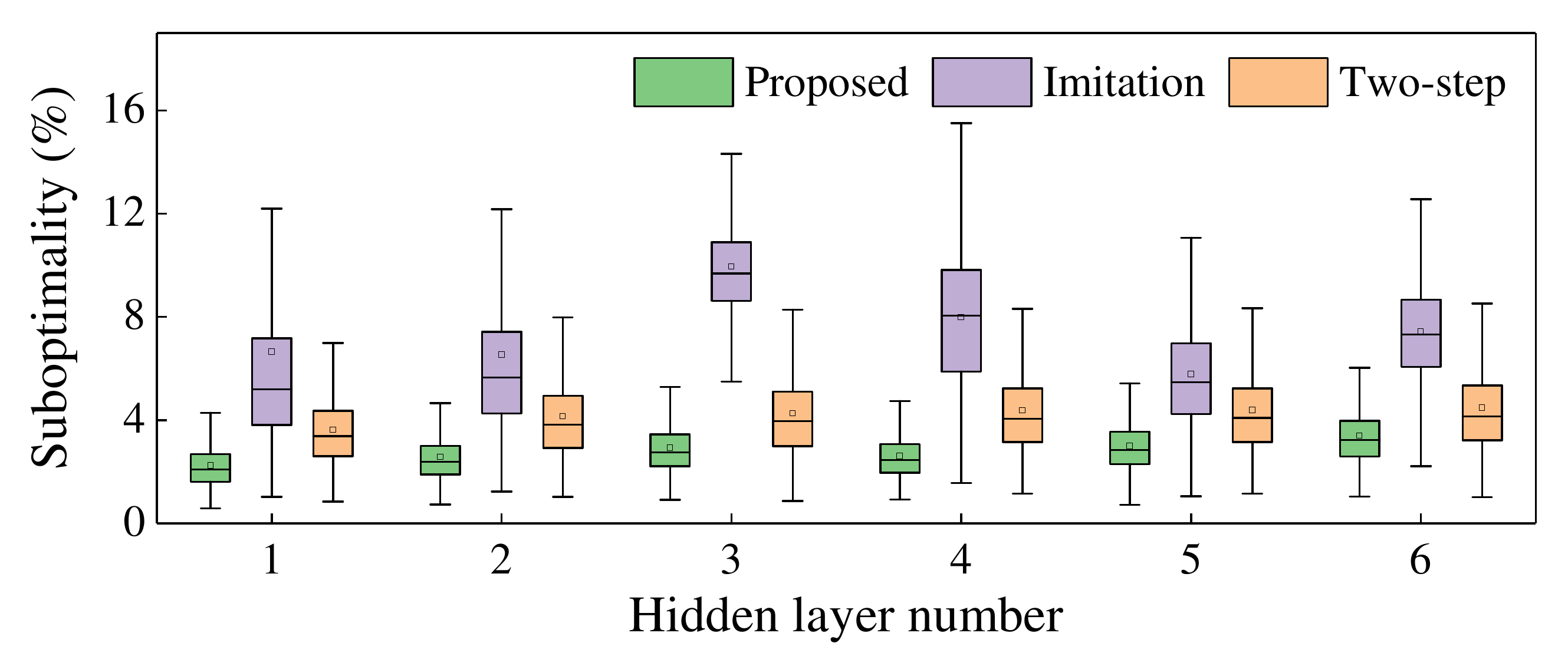}
 	\vspace{-8mm}
 	\caption{Normalized suboptimality of all methods with different hidden layer numbers in the 118-bus test case.}
	\label{fig_regret_layerNum}
	\vspace{-2mm}
\end{figure}

We further investigate the impact of the hidden layer number. While keeping the total neuron number fixed at 60, we vary the hidden layer number from 1 to 6, setting the neurons per layer to 60, 30, 20, 15, 12, and 10, respectively. The suboptimality results for the different methods are shown in Fig.~\ref{fig_regret_layerNum}. For the proposed method, a slight growing trend in suboptimality is observed as the hidden layer number increases. This trend results from the fact that, in our setup, shallower networks have more learnable parameters than deeper ones.\footnote{The total number of learnable parameters for networks with 1, 2, 3, 4, 5, and 6 hidden layers are 14160, 7980, 5520, 4215, 3408, and 2860, respectively.} If the NNs are well-trained, having more learnable parameters (as in the shallower networks) can reduce the generalization error. Thus, shallower networks may achieve lower suboptimality. Moreover, the proposed method directly employs the RLD objective as its training loss, inherently minimizing suboptimality. In contrast, benchmarks \textsf{Imitation} and \textsf{Two-step} use MSE as their training loss, which can not accurately measure suboptimality. Therefore, even when the MSE is minimized, these benchmarks may still exhibit higher suboptimality compared to the proposed method.

\subsubsection{Sample size}
\begin{figure}[h]
	\vspace{-4mm}
	\centering
	{\includegraphics[width=1\columnwidth]{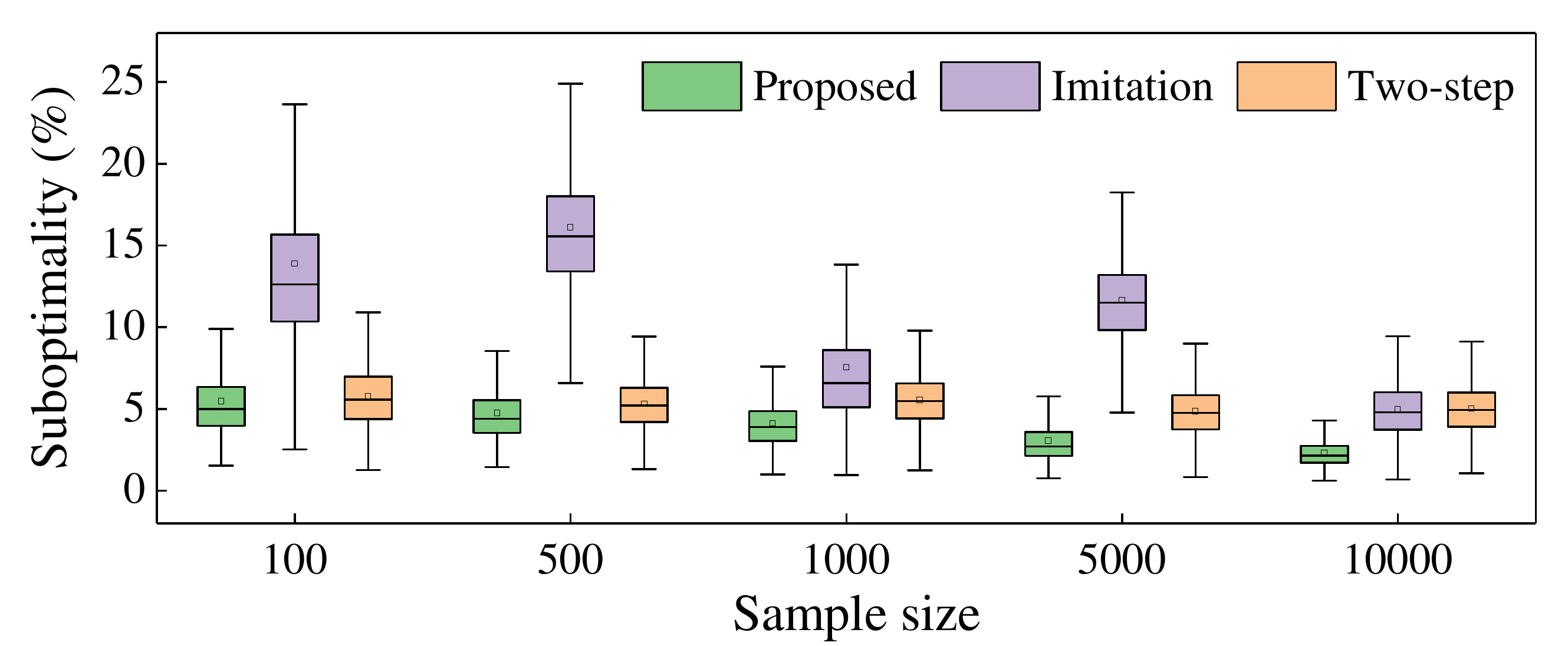}}
	\vspace{-8mm}
 	\caption{Normalized suboptimality of all methods with different training sample sizes in the 118-bus test case.}
	\label{fig_regret_sampleNum}
	\vspace{-4mm}
\end{figure}

Fig. \ref{fig_regret_sampleNum} compares the normalized suboptimality of all methods with different sample sizes. Similar to the results in Fig. \ref{fig_MSE_neuronNum}, the proposed method always outperforms the rest ones by achieving the lowest normalized suboptimality. {Moreover, with an increasing sample size, the suboptimality of the proposed method decreases. According to our proposed generalization performance guarantee, i.e., Theorem \ref{prop_2}, the PAC bound of our method gradually approaches zero with the increase of the sample size.  Other methods do not display a similar trend in suboptimality. Moreover, our proposed method trains the NN by minimizing the RLD loss \eqref{eqn_loss} without requiring labelled data. In other words, the NN is trained in a self-supervised manner. Thus, our approach enables the easy generation of synthetic day-ahead observable information to further enhance the performance of our method.}

\section{Conclusions} \label{sec_conclusion}

This paper presents a learning-based approach to solve the network RLD problem with provably near-optimal performance. It begins by designing a data-driven formulation for the network RLD problem, which aims to learn a decision rule from historical data. This decision rule directly maps day-ahead auxiliary information to cost-effective dispatch decisions for the future delivery interval, so it can bypass the inherent ``mismatch'' issue of the widely adopted predict-then-optimize paradigm. We then propose neural RLD as a novel solution method for this data-driven formulation. It employs a L2-regularized NN to learn the decision rule, thereby transforming the above data-driven formulation as a NN training task with a specialized training loss. The gradient of this loss is explicitly derived, allowing for efficient completion of the training task using SGD.  A generalization performance guarantee is further developed based on uniform convergence, which provides a PAC upper bound for the suboptimality of the neural RLD. Our numerical experiments on IEEE 5-bus, 118-bus, 300-bus, and 1354-bus test systems demonstrate the neural RLD's desirable performance in terms of convergence, suboptimality, and computational efficiency. Simulations also confirm the asymptotic optimality of our method, implying that its suboptimality diminishes to zero with high probability as more samples are utilized.

The proposed neural RLD provides an effective solution to the two-stage network RLD problem. However, its applicability to the original multi-stage stochastic control formulation remains limited. Furthermore, our current theoretical guarantees are restricted to NNs consisting solely of fully connected layers. To address these limitations, future work will focus on extending the neural RLD framework to the multi-stage formulation and developing theoretical guarantees for more advanced learning models. 
\appendices
\setcounter{table}{0}

\section{} \label{app_3}
{\emph{Proof of Proposition \ref{prop_1}}: This proposition is an application of the envelope theorem \cite{milgrom2002envelope}. A similar idea has also been reported in \cite{liang2022operation}.} Specifically, the optimal objective of the dual problem of \eqref{eqn_RLD_2nd} is expressed as:
\begin{align}
V(\mathbf u, \mathbf d) = &\left(\bm \lambda^\star \right)^\intercal \mathbf u - (\bm \lambda^\star)^\intercal \mathbf d - \left(\bm \nu_{-}^\star + \bm \nu_{+}^\star\right)^\intercal \mathbf f^\mathrm{max}.
\end{align}
Since \eqref{eqn_RLD_2nd} can be reformulated as a linear program, the strong duality holds, so
\begin{align}
V(\mathbf u, \mathbf d) = Q(\mathbf u, \mathbf d).
\end{align}
Then, the derivative $\nabla_{\mathbf{u}} Q$ can be expressed as (\ref{eqn_gradient}).

\section{} \label{app_1}

\begin{figure}[h]

    \vspace{-2mm}
	\centering
	{\includegraphics[width=0.9\columnwidth]{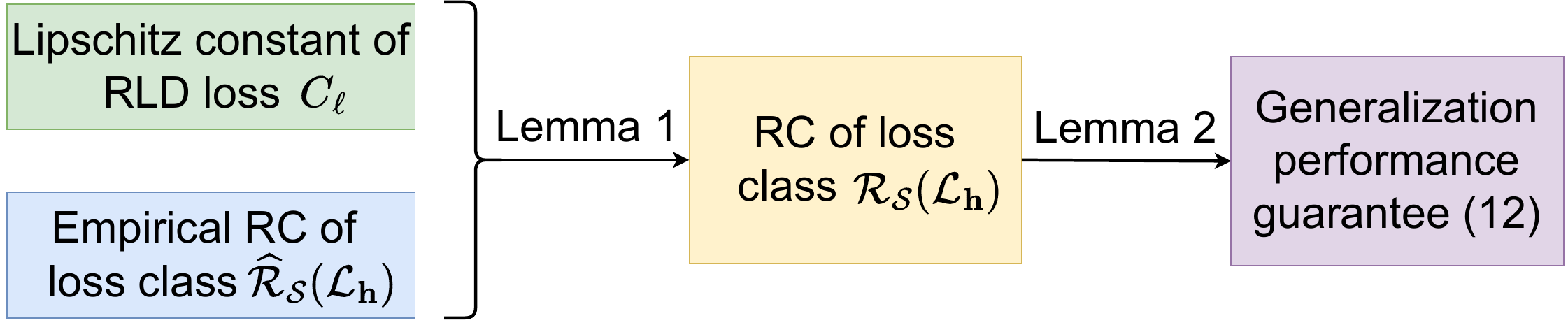}}
	\vspace{-2mm}
 	\caption{Detailed steps for proving Theorem \ref{prop_2}, where RC represents Rademacher complexity.}
	\label{fig_proof}
	\vspace{-2mm}
\end{figure}

\emph{Proof of Theorem \ref{prop_2}}: Fig. \ref{fig_proof} outlines the detailed steps for proving Theorem \ref{prop_2}, as follows: i) Prove the Lipschitz continuity of the RLD loss \eqref{eqn_loss} and calculate the corresponding Lipschitz constant, ii) Estimate the empirical and true Rademacher complexities of the loss class, and iii) Derive the PAC upper bound as stated in Theorem \ref{prop_2}.

\subsubsection{Lipschitz continuity of the RLD loss}

The solution function of problem \eqref{eqn_RLD_2nd} maps the first-stage decision \( \mathbf u \) and net demand \( \mathbf d \) to the second-stage optimal solution $\mathbf g^\star$. According to \cite[Theorem 2.2]{mangasarian1987lipschitz}, this solution function is Lipschitz continuous, and the following inequality holds:
\begin{align}
\Vert (\mathbf{g}^\star)^{(1)} - (\mathbf{g}^\star)^{(2)} \Vert_\infty \leq C_{\mathbf{g}} \Vert \mathbf{u}^{(1)} - \mathbf{u}^{(2)} \Vert_2, 
\end{align}
where $\mathbf{u}^{(1)}$ and $\mathbf{u}^{(2)}$ are two different first-stage decisions, and $\mathbf{g}^\star)^{(1)}$ and $(\mathbf{g}^\star)^{(2)}$ are the corresponding optimal second-stage decisions. Parameter $C_{\mathbf{g}}$ is the Lipschitz constant of this solution function, and its detailed formulation is provided in  \cite[Theorem 2.2]{mangasarian1987lipschitz}.  
Then, we have:
\begin{align}
&\left|Q(\mathbf{u}^{(1)}, \mathbf d^{(1)}) - Q(\mathbf{u}^{(2)}, \mathbf d^{(2)})\right| \notag \\
= &\left|\bm \beta^\intercal\left((\mathbf{g}^\star)^{(1)}\right)_+ - \bm \beta^\intercal \left((\mathbf{g}^\star)^{(2)}\right)_+ \right| \notag \\
\leq & \Vert \bm \beta \Vert_2 \cdot \Vert (\mathbf{g}^\star)^{(1)} - (\mathbf{g}^\star)^{(2)} \Vert_2 \notag \\
\leq & \sqrt{\Nbus} \cdot \Vert \bm \beta \Vert_2 \cdot \Vert (\mathbf{g}^\star)^{(1)} - (\mathbf{g}^\star)^{(2)} \Vert_\infty \notag \\
\leq & \sqrt{\Nbus} C_{\mathbf{g}} \cdot  \Vert \bm \beta \Vert_2 \cdot \Vert \mathbf{u}^{(1)} - \mathbf{u}^{(2)} \Vert_2, 
\end{align}
where the first inequality holds because function $\mathbf{g} \rightarrow \bm \beta^\intercal(\mathbf{g})_+$ is $\Vert \bm \beta \Vert_2$-Lipschitz continuous; the second inequality holds by the relation between infinite and L2 norms.
Finally, the Lipschitz continuity of $\ell(\mathbf u, \mathbf d)$ can be proved by:
\begin{align}
&\left|\ell(\mathbf{u}^{(1)}, \mathbf d^{(1)}) - \ell(\mathbf{u}^{(2)}, \mathbf d^{(2)})\right| \notag \\ 
= & \left| \bm \alpha^\intercal \mathbf{u}^{(1)} + Q(\mathbf{u}^{(1)}, \mathbf d^{(1)}) - \bm \alpha^\intercal \mathbf{u}^{(2)} - Q(\mathbf{u}^{(2)}, \mathbf d^{(2)}) \right| \notag \\
\leq & \underbrace{\left(\Vert \bm \alpha \Vert_2 + \sqrt{\Nbus} C_{\mathbf{g}} \Vert \bm \beta \Vert_2 \right)}_{C_{\ell}}  \Vert \mathbf{u}^{(1)} - \mathbf{u}^{(2)} \Vert_2. \label{eqn_l2u}
\end{align}
Note \eqref{eqn_l2u} also gives the definition of $C_{\ell}$ in Theorem \ref{prop_2}. When the L2-regularized NN $\mathbf h$ is used to predict the dispatch decision $\mathbf{u}$, we have:
\begin{align}
\Vert \mathbf{u}^{(1)} - \mathbf{u}^{(2)} \Vert_2 &= \Vert \mathbf h(\mathbf{x}^{(1)}) - \mathbf h(\mathbf{x}^{(2)}) \Vert_2. \label{eqn_u2h}
\end{align}
By substituting (\ref{eqn_u2h}) into (\ref{eqn_l2u}), we obtain:
\begin{align}
&\left|\ell(\mathbf{u}^{(1)}, \mathbf d^{(1)}) - \ell(\mathbf{u}^{(2)}, \mathbf d^{(2)})\right|\leq C_{\ell} \Vert \mathbf h(\mathbf{x}^{(1)}) - \mathbf h(\mathbf{x}^{(2)}) \Vert_2. \label{eqn_l2h}
\end{align}
The above inequality shows that the mapping $\psi(\cdot): \mathbf h(\mathbf x) \mapsto \ell(\mathbf h(\mathbf x), \mathbf d)$ is $C_{\ell}$-Lipschitz continuous. Then, according to

\subsubsection{Upper bound of Rademacher complexity}
We first define the following new hypothesis class for convenience:
\begin{align}
\hspace{-2mm}
\mathcal{H}^k = \left\{
\mathbf h^k(\mathbf x) \left|
\begin{aligned}
&\mathbf h^k(\mathbf x) = \mathbf W^k \phi \left(\cdots \mathbf W^{2} \phi\left(\mathbf W^{1} \mathbf x\right)\right), \\
&\bm \Vert \mathbf w^{k'}_{j} \Vert_2 \leq W^{\max}, \forall j \in \mathcal{J}^{k'}, \\
&\quad \quad \quad \quad \quad \quad \ \ \ \   \forall k' \in [1,\cdots,k],
\end{aligned} \right. \label{eqn_class_layer}
\right\},
\end{align}
where the hypothesis $\mathbf h^k(\mathbf x): \mathbb{R}^{p} \mapsto \mathbb{R}^{J^k}$ represents a function maps $\mathbf{x}$ to the output of the $k$-th hidden layer. 

The empirical Rademacher complexity of the loss class on dataset $\mathcal{S}$ is defined as follows:
\begin{align}
\widehat{\mathcal{R}}_{\mathcal{S}}(\mathcal{L}_{\mathbf{h}}) = \mathbb{E}_{\xi} \left[ \sup_{\mathbf{\ell} \in \mathcal{L}_{\mathbf{h}}} \sum_{m \in \mathcal{M}} \xi^{(m)} \ell(\mathbf h(\mathbf x^{(m)}), \mathbf d^{(m)}) \right], \label{RC_empirical}	
\end{align}
where $\widehat{\mathcal{R}}_{\mathcal{S}}(\mathcal{L}_{\mathbf{h}})$ is calculated based on observed realizations, which differs from the Rademacher complexity $\mathcal{R}_{\mathcal{S}}(\mathcal{L}_{\mathbf h})$. To obtain the upper bound of $\widehat{\mathcal{R}}_{\mathcal{S}}(\mathcal{L}_{\mathbf{h}})$, we first introduce the \emph{vector contraction inequality} in the following lemma:
\begin{lemma} [Vector contraction inequality \cite{mangasarian1987lipschitz}] \label{lemma_1}
    Let $\mathcal{H}$ be a class of vector-valued functions $\mathbf{h}(\mathbf{x}): \mathbb{R}^{p} \mapsto \mathbb{R}^{\Nbus}$ and let $\psi(\cdot): \mathbb{R}^{\Nbus} \mapsto \mathbb{R}$ is a $C_{\ell}$-Lipschitz continuous function. Then, the following vector contraction inequality holds:
\begin{align}
&\mathbb{E}_{\xi} \left( \sup_{\mathbf{h} \in \mathcal{H}} \sum_{m \in \mathcal{M}} \xi^{(m)} \psi(\mathbf h(\mathbf{x}^{(m)})) \right) \notag \\
& \quad \  \leq \sqrt{2} C_{\ell} \mathbb{E}_{\xi} \left( \sup_{h_j \in \mathcal{H}_j} \sum_{m \in \mathcal{M}}\sum_{j =1}^{\Nbus}  \xi^{(m,j)} h_j(\mathbf{x}^{(m)}) \right),
\end{align}
where $\xi^{(m,j)}$ is a random variable and drawn i.i.d.. from a uniform distribution over $\{-1, +1\}$; $h_j(\mathbf{x}):\mathbb{R}^{p} \mapsto \mathbb{R}$ is a mapping from $\mathbf x$ to the $j$-th component of $\mathbf h$; $\mathcal{H}_j$ represents the hypothesis class of $h_j(\mathbf{x})$ defined such that $\mathbf h_j \in \mathcal H_j$.
\end{lemma}

By substituting Lemma \ref{lemma_1} into \eqref{RC_empirical}, we have:
\begin{align}
\hspace{-2mm}
\widehat{\mathcal{R}}_{\mathcal{S}}(\mathcal{L}_{\mathbf{h}}) \leq \sqrt{2} C_{\ell} \mathbb{E}_{\xi}\left[ \sup_{{h}_j \in \mathcal{H}_j} \sum_{m \in \mathcal{M}}\sum_{j =1}^{\Nbus}  \xi^{(m,j)} {h}_j(\mathbf{x}^{(m)}) \right]. \label{eqn_vector_R_1}
\end{align}
Reference \cite{mangasarian1987lipschitz} further proved the following inequality:
\begin{align}
&\mathbb{E}_{\xi}\left[ \sup_{{h}_j \in \mathcal{H}_j} \sum_{m \in \mathcal{M}}\sum_{j =1}^{\Nbus}  \xi^{(m,j)} {h}_j(\mathbf{x}^{(m)}) \right] \notag \\ 
&\leq \sum_{j =1}^{\Nbus} \mathbb{E}_{\xi}\left[ \sup_{{h}_j \in \mathcal{H}_j} \sum_{m \in \mathcal{M}}  \xi^{(m,j)} {h}_j(\mathbf{x}^{(m)}) \right]. \label{eqn_vector_R_2}
\end{align}
Then, this empirical Rademacher complexity is bounded by:
\begin{align}
 \widehat{\mathcal{R}}_{\mathcal{S}}(\mathcal{L}_{\mathbf{h}}) &\leq \sqrt{2} C_{\ell} \sum_{j =1}^{\Nbus} \mathbb{E}_{\xi}\left[ \sup_{{h}_j \in \mathcal{H}_j} \sum_{m \in \mathcal{M}}  \xi^{(m,j)} {h}_j(\mathbf{x}^{(m)}) \right], \label{eqn_UB0}
\end{align}
where $\mathcal{X} = \{\mathbf{x}_i, \forall m \in \mathcal{M}\}$. Observing that the empirical Rademacher complexity of the hypothesis class $\mathcal{H}_j$ on dataset $\mathcal{X}$ is defined as:
\begin{align}
\widehat{\mathcal{R}}_{\mathcal{X}}(\mathcal{H}_j) = \mathbb{E}_{\xi}\left[ \sup_{{h}_j \in \mathcal{H}_j} \sum_{m \in \mathcal{M}}  \xi^{(m,j)} {h}_j(\mathbf{x}^{(m)}) \right].
\end{align}
Eq. \eqref{eqn_UB0} can be further expressed as:
\begin{align}
 \widehat{\mathcal{R}}_{\mathcal{S}}(\mathcal{L}_{\mathbf{h}}) \leq \sqrt{2} C_{\ell} \sum_{j =1}^{\Nbus} \widehat{\mathcal{R}}_{\mathcal{X}}(\mathcal{H}_j). \label{eqn_UB1}
\end{align} 

Before calculating $\widehat{\mathcal{R}}_{\mathcal{X}}(\mathcal{H}_j)$, we first introduce the positive homogeneity of ReLU $\phi(\cdot)$ as follows:
\begin{align}
a \cdot \phi(\mathbf x) =  \phi(a  \cdot \mathbf x), \ \forall a \geq 0. \label{eqn_positive homogeneity}
\end{align}
Meanwhile, from the L2 regularization, we have:
\begin{align}
\Vert \mathbf w^{k}_{j} \Vert_2 \leq W^\mathrm{max} \quad \Rightarrow \quad \Vert \mathbf w^{k}_{j} \Vert_\infty \leq  \sqrt{W^\mathrm{max}}. \label{eqn_regularization}
\end{align}
We further define a new variable $\overline{\mathbf w}^{k}_{j} = \frac{\mathbf w^{k}_{j}}{{\Vert \mathbf w^{k}_{j} \Vert_2}}$. Obviously, we have $\Vert \overline{\mathbf w}^{k}_{j} \Vert_2 = 1$. Then, according to the definition of $h_j(\mathbf x)$, the empirical Rademacher complexity $\widehat{\mathcal{R}}_{\mathcal{X}}(\mathcal{H}_j)$ can be expressed as:
\begin{align}
&\widehat{\mathcal{R}}_{\mathcal{X}}(\mathcal{H}_j) = \frac{1}{M} \mathbb{E}_{\xi} \left[ \sup_{ \mathbf{h}^{K-2} \in \mathcal{H}^{K-2}}  \sum_{m \in \mathcal{M}} \xi^{(m)}  \left( \sum_{j' \in \mathcal{J}^{K-1}} w^{K}_{j,j'} \right.\right. \notag \\
& \quad \quad \quad \quad \quad \quad \quad \quad \quad \left. \left. \cdot \phi\left((\mathbf w^{K-1}_{j'})^\intercal \phi(\mathbf h^{K-2}(\mathbf{x}^{(m)}))\right) \right)\right] \notag \\
&= \frac{1}{M} \mathbb{E}_{\xi} \left[ \sup_{ \mathbf h^{K-2} \in \mathcal{H}^{K-2}} \sum_{j' \in \mathcal{J}^{K-1}} \underbrace{w^{K}_{j,j'} \left\Vert \mathbf w^{K-1}_{j'} \right\Vert_2}_{a_{j'}} \right. \notag \\
& \quad \quad \quad \quad \quad \quad    \left. \cdot \underbrace{\sum_{m \in \mathcal{M}} \xi^{(m)} \phi\left( (\overline{\mathbf w}^{K-1}_{j'})^\intercal  \phi(\mathbf h^{(K-2)})\right)}_{b_{j'}} \right], \label{eqn_ERC_1}
\end{align}
where $w^{K}_{j,j'}$ is the $j'$ component of vector $\mathbf w^{K}_{j}$. The second inequality holds due to the homogeneity of ReLU. Since $\sum_{j'} a_{j'} b_{j'} \leq \sum_{j'} |a_{j'}| \max_{j''}|{b_{{j''}}}|$ holds for any $a_{j'}$ and $b_{j'}$, the value of $\widehat{\mathcal{R}}_{\mathcal{X}}(\mathcal{H}_j)$ is bounded by:
\begin{align}
&\widehat{\mathcal{R}}_{\mathcal{X}}(\mathcal{H}_j) \leq \frac{1}{M} \mathbb{E}_{\xi} \left[ \sup_{ \mathbf h^{K-2} \in \mathcal{H}^{K-2}} \sum_{j' \in \mathcal{J}^{K-1}} \underbrace{\left| w^{K}_{j,j'}\right| \left\Vert \mathbf w^{K-1}_{j'} \right\Vert_2}_{|a_{j'}|} \right. \notag \\
& \quad \quad \quad \quad  \left. \cdot \underbrace{\max_{j''} \left|\sum_{m \in \mathcal{M}} \xi^{(m)} \phi\left( (\overline{\mathbf w}^{K-1}_{j''})^\intercal  \phi(\mathbf h^{K-2})\right)\right|}_{\max_{j''}|{b_{j''}}|} \right].
\end{align}
According to \eqref{eqn_regularization}, we further have:
\begin{align}
&\widehat{\mathcal{R}}_{\mathcal{X}}(\mathcal{H}_j) \notag \\
 &\leq  \frac{(W^\mathrm{max})^{\frac{3}{2}}}{M} \mathbb{E}_{\xi} \left[ \sup_{\substack{ \mathbf h^{K-2} \in \mathcal{H}^{K-2},\\ \left\Vert \overline{\mathbf w} \right\Vert_2=1,}} \left|\sum_{m \in \mathcal{M}} \xi^{(m)} \phi\left( \overline{\mathbf w}^\intercal  \phi(\mathbf h^{K-2})\right)\right| \right] \notag \\
&\leq \frac{(W^\mathrm{max})^{\frac{3}{2}}}{M} \mathbb{E}_{\xi} \left[ \sup_{\substack{ \mathbf h^{K-2} \in \mathcal{H}^{K-2},\\ \left\Vert \overline{\mathbf w} \right\Vert_2 \leq 1,}} \left|\sum_{m \in \mathcal{M}} \xi^{(m)} \phi\left( \overline{\mathbf w}^\intercal  \phi(\mathbf h^{K-2})\right)\right| \right]. \label{eqn_ERC_2}
\end{align}
In the first inequality above, we replace $\overline{\mathbf w}^{K-1}_{j'}$ by a new variable $\overline{\mathbf w}$, which has the same dimension as $\mathbf w^{K-1}_{j'}$.
this replacement is equivalent because $\overline{\mathbf w}^{K-1}_{j'}$ and $\overline{\mathbf w}$ are only constrained by $\left\Vert \overline{\mathbf w} \right\Vert_2 = 1$ and $\Vert \mathbf w^{K-1}_{j'}\Vert_2 = 1$, respectively. The second inequality holds because it enlarges the hypothesis class. Based on [Lemma 5.12] in \cite{ma2022lecture}, we can remove the absolute operator in \eqref{eqn_ERC_2} and obtain:
\begin{align}
&\widehat{\mathcal{R}}_{\mathcal{X}}(\mathcal{H}_j) \notag \\
& \leq \frac{2(W^\mathrm{max})^{\frac{3}{2}}}{M} \mathbb{E}_{\xi} \left[ \sup_{\substack{ \mathbf h^{K-2} \in \mathcal{H}^{K-2},\\ \left\Vert \overline{\mathbf w} \right\Vert_2 \leq 1,}} \sum_{m \in \mathcal{M}} \xi^{(m)} \phi\left( \overline{\mathbf w}^\intercal  \phi(\mathbf h^{K-2})\right) \right] \notag \\
& \leq \frac{2(W^\mathrm{max})^{\frac{3}{2}}}{M} \mathbb{E}_{\xi} \left[ \sup_{\substack{ \mathbf h^{K-2} \in \mathcal{H}^{K-2},\\ \left\Vert \overline{\mathbf w} \right\Vert_2 \leq 1,}} \sum_{m \in \mathcal{M}} \xi^{(m)} \left( \overline{\mathbf w}^\intercal  \phi(\mathbf h^{K-2})\right) \right] \notag \\
& = 2(W^\mathrm{max})^{\frac{3}{2}} \widehat{\mathcal{R}}_{\mathcal{X}}(\overline{\mathcal{H}}^{K-1}), \label{eqn_proof}
\end{align}
where the definition of class $\overline{\mathcal{H}}^{k}$ is almost the same as ${\mathcal{H}}^{k}$ defined in \eqref{eqn_class_layer} except for replacing the L2 regulation constraint $\Vert \mathbf w^{(k)}_j\Vert_2 \leq W^\mathrm{max}$ with $\left\Vert \overline{\mathbf w} \right\Vert_2 \leq 1$. The second inequality holds because $\phi(\cdot)$ is 1-Lipschitz continuous. The last equality holds due to the definition of the empirical Rademacher complexity. In a similar way, we can also obtain: 
\begin{align}
\widehat{\mathcal{R}}_{\mathcal{X}}(\overline{\mathcal{H}}^{k}) \leq 2W^\mathrm{max} \widehat{\mathcal{R}}_{\mathcal{X}}(\overline{\mathcal{H}}^{k-1}), \forall k \in \{2,\cdots, K-1\}, \label{eqn_proof_2}
\end{align}
where the coefficient is $2W^\mathrm{max}$ instead of $2(W^\mathrm{max})^{\frac{3}{2}}$ because the infinity norm of $\overline{\mathbf w}$ is bounded by one instead of $\sqrt{W^\mathrm{max}}$. 
By substituting (\ref{eqn_proof}) and (\ref{eqn_proof_2}) recursively, we obtain:
\begin{align}
\widehat{\mathcal{R}}_{\mathcal{X}}(\mathcal{H}_j) &\leq \sqrt{W^\mathrm{max}}\cdot 2W^\mathrm{max} \cdot  \widehat{\mathcal{R}}_{\mathcal{X}}(\overline{\mathcal{H}}^{K-1}) \notag \\
&\leq \sqrt{W^\mathrm{max}} \cdot (2W^\mathrm{max})^{2}  \cdot   \widehat{\mathcal{R}}_{\mathcal{X}}(\overline{\mathcal{H}}^{K-2}) \notag \\
& \cdots \notag \\
& \leq \sqrt{W^\mathrm{max}}\cdot (2W^\mathrm{max})^{K-1}  \cdot   \widehat{\mathcal{R}}_{\mathcal{X}}(\overline{\mathcal{H}}^{1}). \label{eqn_NN2linear}
\end{align}

The Rademacher complexity of the loss is the expectation of the empirical one over the probability distribution of $\mathbf x$. By substituting (\ref{eqn_NN2linear}) into \eqref{eqn_UB1}, we can bound this Rademacher complexity via:
\begin{align}
\mathcal{R}_{\mathcal{S}}(\mathcal{L}_{\mathbf h}) 
 &= \mathbb{E}_{\mathbf x}\left(\widehat{\mathcal{R}}_{\mathcal{S}}(\mathcal{L}_{\mathbf{h}})\right) \notag \\ 
 &\leq \sqrt{2} C_{\ell} \cdot \mathbb{E}_{\mathbf x}\left(\widehat{\mathcal{R}}_{\mathcal{X}}(\mathcal{H}_j^{K})\right) \notag \\
&\leq \sqrt{2} C_{\ell} \cdot \sqrt{W^\mathrm{max}}\cdot (2W^\mathrm{max})^{K-1}  \mathbb{E}_{\mathbf x}\left(\widehat{\mathcal{R}}_{\mathcal{X}}(\overline{\mathcal{H}}^{1})\right) \notag \\
&\leq \sqrt{2} C_{\ell} \cdot \sqrt{W^\mathrm{max}}\cdot (2W^\mathrm{max})^{K-1}   \mathcal{R}_\mathcal{S}(\overline{\mathcal{H}}^{1,L2}), \label{eqn_approximation_3}
\end{align}
where $\overline{\mathcal{H}}^{1,L2} = \{\mathbf x \rightarrow \overline{\mathbf w}^\intercal \mathbf x, \Vert \overline{\mathbf w} \Vert_2 \leq 1 \}$. According to \cite[Theorem 3]{kakade2008complexity}, its Rademacher complexity is bounded by:
\begin{align}
\mathcal{R}_\mathcal{S}(\overline{\mathcal{H}}^{1,L2}) \leq \frac{X^\mathrm{max}}{\sqrt{M}}. \label{eqn_approximation_4}
\end{align}
By substituting (\ref{eqn_approximation_4}) into (\ref{eqn_approximation_3}), we obtain an upper bound of the Rademacher complexity $\mathcal{R}_{\mathcal{S}}(\mathcal{L}_{\mathbf h})$, as follows:
\begin{align}
\mathcal{R}_{\mathcal{S}}(\mathcal{L}_{\mathbf h}) \leq \frac{ (2W^\mathrm{max})^{K-\frac{1}{2}} C_{\ell} X^\mathrm{max}}{\sqrt{M}}. \label{eqn_R_bound2}
\end{align}

\subsubsection{PAC upper bound of the excess cost}

The following Lemma gives a PAC bound for the excess cost of the hypothesis $\widehat{\mathbf h}$ based on uniform
convergence.
\begin{lemma} [Generalization bound via uniform
convergence \cite{bartlett2002rademacher}] 
When Assumptions \textbf{A1} and \textbf{A2} hold, with probability at least $1-\delta$ for any small $\delta \in (0,1)$, the excess cost of hypothesis $\widehat{\mathbf h}$ is bounded by:
\begin{align}
\Delta L(\widehat{\mathbf h}) \leq 4 \mathcal{R}_{\mathcal{S}}(\mathcal{L}_{\mathbf h}) + \sqrt{\frac{2\ln(2/\delta)}{M}}. \label{eqn_EC_bound}
\end{align} 
\end{lemma}
By substituting \eqref{eqn_R_bound2} into \eqref{eqn_EC_bound}, Theorem \ref{prop_2} can be proved.

\cg{\subsubsection{Extension to other types of NNs}
The performance guarantee established for L2-regularized NNs can be extended to other types of NNs with minor modifications to the derivation process. Taking L1-regularized NNs as an example, the L1-regularization constraint implies:
\begin{align}
\|\mathbf{w}^{k}_{j}\|_1 \leq W^\mathrm{max} \quad \Rightarrow \quad \|\mathbf{w}^{k}_{j}\|_\infty \leq \|\mathbf{w}^{k}_{j}\|_1 \leq W^\mathrm{max}.
\end{align}
Following Eqs. \eqref{eqn_ERC_1}-\eqref{eqn_ERC_2}, the empirical Rademacher complexity \(\widehat{\mathcal{R}}_{\mathcal{X}}(\mathcal{H}_j)\) is bounded by:
\begin{align}
&\widehat{\mathcal{R}}_{\mathcal{X}}(\mathcal{H}_j) \notag \\ 
&\leq  \frac{(W^\mathrm{max})^{2}}{M} \mathbb{E}_{\xi} \left[ \sup_{\substack{\mathbf{h}^{K-2} \in \mathcal{H}^{K-2},\\ \|\overline{\mathbf{w}}\|_2 = 1}} \left|\sum_{m \in \mathcal{M}} \xi^{(m)} \phi\left( \overline{\mathbf{w}}^\intercal \phi(\mathbf{h}^{K-2})\right)\right| \right].
\end{align}
Further applying Eqs. \eqref{eqn_proof}--\eqref{eqn_approximation_3}, we derive:
\begin{align}
\mathcal{R}_{\mathcal{S}}(\mathcal{L}_{\mathbf{h}}) 
\leq \sqrt{2} C_{\ell} \cdot (2W^\mathrm{max})^{K} \mathcal{R}_{\mathcal{S}}(\overline{\mathcal{H}}^{1,L1}),
\end{align}
where the hypothesis class is defined as \(\overline{\mathcal{H}}^{1,L1} = \{\mathbf{x} \mapsto \overline{\mathbf{w}}^\intercal \mathbf{x} : \|\overline{\mathbf{w}}\|_1 \leq 1 \}\). According to \cite[Theorem 3]{kakade2008complexity}, its Rademacher complexity is bounded by:
\begin{align}
\mathcal{R}_{\mathcal{S}}(\overline{\mathcal{H}}^{1,L1}) \leq X^\mathrm{max} \sqrt{\frac{\ln p}{M}},
\end{align}
where \(p\) is the dimensionality of the observable feature vector \(\mathbf{x}\). The remaining steps in the derivation directly follow the process used for L2-regularized NNs to establish the performance guarantee for L1-regularized NNs.
}

\cg{References \cite{joukovsky2021generalization,valleperez2020generalizationboundsdeeplearning,truong2024rademachercomplexitybasedgeneralizationbounds} provide bounds for the empirical Rademacher complexity \(\widehat{\mathcal{R}}_{\mathcal{X}}(\mathcal{H}_j)\) when RNNs or CNNs are used as the hypothesis classes for prediction. Substituting these bounds into \eqref{eqn_UB1} yields the Rademacher complexity bound for the loss class \(\mathcal{R}_{\mathcal{S}}(\mathcal{L}_{\mathbf{h}})\). By further substituting this bound into \eqref{eqn_EC_bound}, we can obtain the extension of the performance guarantee to RNNs and CNNs.}

\section{} \label{app_4}

\subsubsection{Extension to the AC formulation}
After replacing the DC power flow model in \eqref{eqn_DCOPF}-\eqref{eqn_capacity} with the AC formulation, the second-stage problem becomes non-convex. Nevertheless, convex relaxations of the AC power flow model have been proposed \cite{9862539,10508977}. According to the envelope theorem mentioned in Proposition \ref{prop_1}, the gradient of an optimization problem's Lagrangian at saddle points with respect to certain parameters can be explicitly expressed using dual variables. If our second-stage problem has a tight convex relaxation, its Lagrangian at the saddle points will equal its optimal objective value. Therefore, the gradient of the RLD loss can be computed using dual variables, enabling effective NN training. With this tight convex relaxation, the value function $Q(\mathbf{u}, \mathbf{d})$ remains Lipschitz continuous, as shown in \cite[Theorem 5.2]{still2018lectures}. Thus, the RLD loss also retains its Lipschitz continuity. By replacing the constant $C_\ell$ in \eqref{eqn_EC_bound_final} with the Lipschitz constant of the new $Q(\mathbf{u}, \mathbf{d})$, we can apply the same procedure detailed in Appendix \ref{app_1} to establish the generalization performance guarantee for the AC case.

\subsubsection{Extension to incorporating multiple periods and additional operation constraints}
When ramping constraints, production limits, and generation reduction costs are considered, the RLD problem becomes
\begin{align}
	\min_{\mathbf u_t, \forall t \in \mathcal{T}} \quad \left(\sum_{t \in \mathcal{T}} \bm \alpha_t^\intercal \mathbf u_t\right) + \mathbb{E}_{\mathbf d \sim \mathcal{D}(\mathbf x)}[Q(\mathbf u, \mathbf d)], \label{eqn_RLD_realistic}
\end{align}
where we introduce a temporal dimension with time intervals indexed by $t \in \mathcal{T}$ to account for ramping constraints. The second-stage cost $Q(\mathbf u, \mathbf d)$ represents the real-time dispatch cost given the trajectories of first-stage decision $\mathbf u = [\mathbf u_t, \forall t \in \mathcal{T}]$ and realized net demand $\mathbf d = [\mathbf d_t, \forall t \in \mathcal{T}]$:
\begin{subequations} 
\begin{align} 
	Q(\mathbf u, \mathbf d&) =  \notag \\
	\min_{(\mathbf g_t, \bm \theta_t)_{t \in \mathcal{T}}} \ & \sum_{t \in \mathcal{T}} \left((\bm \beta_t^+)^\intercal \mathbf g_t^+ + (\bm \beta_t^-)^\intercal \mathbf g_t^-\right),  \label{eqn_obj_realistic} \\
		 \mbox{s.t.} \quad  & \mathbf g_t^+ - \mathbf  g_t^- = \mathbf g_t, \mathbf g_t^+ \geq \bm 0, \mathbf g_t^- \geq \bm 0, \ \forall t \in \mathcal{T}, \label{eqn_cost_realistic}\\
		 & \underline{\mathbf u} \leq \mathbf u_t + \mathbf g_t \leq \overline{\mathbf u}, \ \forall t \in \mathcal{T}, \label{eqn_generation_limit_realistic} \\
		 &  -\mathbf r^\mathrm{down} \leq (\mathbf u_t + \mathbf g_t) - (\mathbf u_{t-1} + \mathbf g_{t-1}) \notag \\
		  & \quad \quad\quad\quad\quad\quad\quad\quad\quad\quad\quad \leq \mathbf r^\mathrm{up}, \ \forall t \in \mathcal{T}, \label{eqn_ramping_realistic} \\
		 & \mathbf u_t + \mathbf g_t - \mathbf d_t = \mathbf B \bm \theta_t, \ \forall t \in \mathcal{T}, \label{eqn_DCOPF_realistic} \\
  & -\mathbf f^\mathrm{max} \leq \mathbf F \bm \theta_t \leq \mathbf f^\mathrm{max},  \ \forall t \in \mathcal{T}, \label{eqn_capacity_realistic} 
\end{align}	
\end{subequations}
where $\mathbf g_t^+$ and $\mathbf g_t^-$ represent real-time generation adjustments, with costs captured by $\bm \beta_t^+$ and $\bm \beta_t^-$ for each time period $t$. Constraints \eqref{eqn_generation_limit_realistic} and \eqref{eqn_ramping_realistic} enforce production and ramping limits, while \eqref{eqn_DCOPF_realistic} and \eqref{eqn_capacity_realistic} describe power flow and branch capacity limits. All operational constraints are imposed in the second stage, ensuring that real-time generation adjustments are feasible and implementable. 

Since the second-stage problem remains linear, the gradient $\nabla_{\bm \phi} Q_t$ can be explicitly expressed using dual variables according to Proposition \ref{prop_1}, allowing for effective training. The generalization performance guarantee can also be derived in the same manner as outlined in Appendix \ref{app_1}, as the linear structure of the second-stage problem remains unchanged.

\bibliographystyle{IEEEtran}
\bibliography{ref}

\begin{thebibliography}{10}
\providecommand{\url}[1]{#1}
\csname url@samestyle\endcsname
\providecommand{\newblock}{\relax}
\providecommand{\bibinfo}[2]{#2}
\providecommand{\BIBentrySTDinterwordspacing}{\spaceskip=0pt\relax}
\providecommand{\BIBentryALTinterwordstretchfactor}{4}
\providecommand{\BIBentryALTinterwordspacing}{\spaceskip=\fontdimen2\font plus
\BIBentryALTinterwordstretchfactor\fontdimen3\font minus
  \fontdimen4\font\relax}
\providecommand{\BIBforeignlanguage}[2]{{%
\expandafter\ifx\csname l@#1\endcsname\relax
\typeout{** WARNING: IEEEtran.bst: No hyphenation pattern has been}%
\typeout{** loaded for the language `#1'. Using the pattern for}%
\typeout{** the default language instead.}%
\else
\language=\csname l@#1\endcsname
\fi
#2}}
\providecommand{\BIBdecl}{\relax}
\BIBdecl

\bibitem{javed2023impact}
M.~S. Javed, J.~Jurasz, M.~Guezgouz, F.~A. Canales, T.~H. Ruggles, and T.~Ma,
  ``Impact of multi-annual renewable energy variability on the optimal sizing
  of off-grid systems,'' \emph{Renewable and Sustainable Energy Reviews}, vol.
  183, p. 113514, 2023.

\bibitem{zhao2022sustainable}
N.~Zhao and F.~You, ``Sustainable power systems operations under renewable
  energy induced disjunctive uncertainties via machine learning-based robust
  optimization,'' \emph{Renewable and Sustainable Energy Reviews}, vol. 161, p.
  112428, 2022.

\bibitem{zakaria2020uncertainty}
A.~Zakaria, F.~B. Ismail, M.~H. Lipu, and M.~A. Hannan, ``Uncertainty models
  for stochastic optimization in renewable energy applications,''
  \emph{Renewable Energy}, vol. 145, pp. 1543--1571, 2020.

\bibitem{5618534}
P.~P. Varaiya, F.~F. Wu, and J.~W. Bialek, ``Smart operation of smart grid:
  Risk-limiting dispatch,'' \emph{Proceedings of the IEEE}, vol.~99, no.~1, pp.
  40--57, 2011.

\bibitem{rajagopal2013risk}
R.~Rajagopal, E.~Bitar, P.~Varaiya, and F.~Wu, ``Risk-limiting dispatch for
  integrating renewable power,'' \emph{International Journal of Electrical
  Power \& Energy Systems}, vol.~44, no.~1, pp. 615--628, 2013.

\bibitem{qin2013risk}
J.~Qin, B.~Zhang, and R.~Rajagopal, ``Risk limiting dispatch with ramping
  constraints,'' in \emph{2013 IEEE International Conference on Smart Grid
  Communications (SmartGridComm)}.\hskip 1em plus 0.5em minus 0.4em\relax IEEE,
  2013, pp. 791--796.

\bibitem{6580485}
J.~Qin, H.-I. Su, and R.~Rajagopal, ``Storage in risk limiting dispatch:
  Control and approximation,'' in \emph{2013 American Control Conference},
  2013, pp. 4202--4208.

\bibitem{6818365}
B.~Zhang, R.~Rajagopal, and D.~Tse, ``Network risk limiting dispatch: Optimal
  control and price of uncertainty,'' \emph{IEEE Transactions on Automatic
  Control}, vol.~59, no.~9, pp. 2442--2456, 2014.

\bibitem{8618917}
J.~Qin, K.~Poolla, and P.~Varaiya, ``Direct data-driven methods for risk
  limiting dispatch,'' in \emph{2018 IEEE Conference on Decision and Control
  (CDC)}, 2018, pp. 3994--4001.

\bibitem{9599541}
L.~V. Krannichfeldt, Y.~Wang, T.~Zufferey, and G.~Hug, ``Online ensemble
  approach for probabilistic wind power forecasting,'' \emph{IEEE Transactions
  on Sustainable Energy}, vol.~13, no.~2, pp. 1221--1233, 2022.

\bibitem{8982039}
H.~Zhang, Y.~Liu, J.~Yan, S.~Han, L.~Li, and Q.~Long, ``Improved deep mixture
  density network for regional wind power probabilistic forecasting,''
  \emph{IEEE Transactions on Power Systems}, vol.~35, no.~4, pp. 2549--2560,
  2020.

\bibitem{8999581}
X.~Fu, Q.~Guo, and H.~Sun, ``Statistical machine learning model for stochastic
  optimal planning of distribution networks considering a dynamic correlation
  and dimension reduction,'' \emph{IEEE Transactions on Smart Grid}, vol.~11,
  no.~4, pp. 2904--2917, 2020.

\bibitem{8651522}
J.~Hu and H.~Li, ``A new clustering approach for scenario reduction in
  multi-stochastic variable programming,'' \emph{IEEE Transactions on Power
  Systems}, vol.~34, no.~5, pp. 3813--3825, 2019.

\bibitem{elmachtoub2022smart}
A.~N. Elmachtoub and P.~Grigas, ``Smart “predict, then optimize”,''
  \emph{Management Science}, vol.~68, no.~1, pp. 9--26, 2022.

\bibitem{9755891}
L.~Sang, Y.~Xu, H.~Long, Q.~Hu, and H.~Sun, ``Electricity price prediction for
  energy storage system arbitrage: A decision-focused approach,'' \emph{IEEE
  Transactions on Smart Grid}, vol.~13, no.~4, pp. 2822--2832, 2022.

\bibitem{zhang2021review}
L.~Zhang, J.~Wen, Y.~Li, J.~Chen, Y.~Ye, Y.~Fu, and W.~Livingood, ``A review of
  machine learning in building load prediction,'' \emph{Applied Energy}, vol.
  285, p. 116452, 2021.

\bibitem{10058008}
G.~Chen, H.~Zhang, and Y.~Song, ``Efficient constraint learning for data-driven
  active distribution network operation,'' \emph{IEEE Transactions on Power
  Systems}, vol.~39, no.~1, pp. 1472--1484, 2024.

\bibitem{10700765}
G.~Chen, J.~Qin, and H.~Zhang, ``Model-free self-supervised learning for
  dispatching distributed energy resources,'' \emph{IEEE Transactions on Smart
  Grid}, vol.~16, no.~2, pp. 1287--1300, 2025.

\bibitem{10256159}
W.~Chen, M.~Tanneau, and P.~Van~Hentenryck, ``End-to-end feasible optimization
  proxies for large-scale economic dispatch,'' \emph{IEEE Transactions on Power
  Systems}, vol.~39, no.~2, pp. 4723--4734, 2024.

\bibitem{NIPS2017_3fc2c60b}
\BIBentryALTinterwordspacing
P.~Donti, B.~Amos, and J.~Z. Kolter, ``Task-based end-to-end model learning in
  stochastic optimization,'' in \emph{Advances in Neural Information Processing
  Systems}, I.~Guyon, U.~V. Luxburg, S.~Bengio, H.~Wallach, R.~Fergus,
  S.~Vishwanathan, and R.~Garnett, Eds., vol.~30.\hskip 1em plus 0.5em minus
  0.4em\relax Curran Associates, Inc., 2017. [Online]. Available:
  \url{https://proceedings.neurips.cc/paper_files/paper/2017/file/3fc2c60b5782f641f76bcefc39fb2392-Paper.pdf}
\BIBentrySTDinterwordspacing

\bibitem{pmlr_v202_rychener23a}
\BIBentryALTinterwordspacing
Y.~Rychener, D.~Kuhn, and T.~Sutter, ``End-to-end learning for stochastic
  optimization: A {B}ayesian perspective,'' in \emph{Proceedings of the 40th
  International Conference on Machine Learning}, ser. Proceedings of Machine
  Learning Research, A.~Krause, E.~Brunskill, K.~Cho, B.~Engelhardt, S.~Sabato,
  and J.~Scarlett, Eds., vol. 202.\hskip 1em plus 0.5em minus 0.4em\relax PMLR,
  23--29 Jul 2023, pp. 29\,455--29\,472. [Online]. Available:
  \url{https://proceedings.mlr.press/v202/rychener23a.html}
\BIBentrySTDinterwordspacing

\bibitem{doi:10.1137/S1052623499363220}
A.~J. Kleywegt, A.~Shapiro, and T.~Homem-de Mello, ``The sample average
  approximation method for stochastic discrete optimization,'' \emph{SIAM
  Journal on Optimization}, vol.~12, no.~2, pp. 479--502, 2002.

\bibitem{shapiro2021lectures}
A.~Shapiro, D.~Dentcheva, and A.~Ruszczynski, \emph{Lectures on stochastic
  programming: modeling and theory}.\hskip 1em plus 0.5em minus 0.4em\relax
  SIAM, 2021.

\bibitem{kim2015guide}
S.~Kim, R.~Pasupathy, and S.~G. Henderson, ``A guide to sample average
  approximation,'' \emph{Handbook of simulation optimization}, pp. 207--243,
  2015.

\bibitem{1266587}
T.~Wu, M.~Rothleder, Z.~Alaywan, and A.~Papalexopoulos, ``Pricing energy and
  ancillary services in integrated market systems by an optimal power flow,''
  \emph{IEEE Transactions on Power Systems}, vol.~19, no.~1, pp. 339--347,
  2004.

\bibitem{li2025frequency}
H.~Li, H.~Zhang, J.~Zhang, Q.~Wu, and C.-K. Wong, ``A frequency-secured
  planning method for integrated electricity-heat microgrids with virtual
  inertia suppliers,'' \emph{Applied Energy}, vol. 377, p. 124540, 2025.

\bibitem{kariuki1996evaluation}
K.~Kariuki and R.~N. Allan, ``Evaluation of reliability worth and value of lost
  load,'' \emph{IEE proceedings-Generation, transmission and distribution},
  vol. 143, no.~2, pp. 171--180, 1996.

\bibitem{mas1995microeconomic}
A.~Mas-Colell, M.~Whinston, and J.~Green, \emph{Microeconomic Theory}, 1995.

\bibitem{montanari2022universality}
A.~Montanari and B.~N. Saeed, ``Universality of empirical risk minimization,''
  in \emph{Conference on Learning Theory}.\hskip 1em plus 0.5em minus
  0.4em\relax PMLR, 2022, pp. 4310--4312.

\bibitem{NEURIPS2020_2000f632}
Y.~Lu and J.~Lu, ``A universal approximation theorem of deep neural networks
  for expressing probability distributions,'' in \emph{Advances in Neural
  Information Processing Systems}, H.~Larochelle, M.~Ranzato, R.~Hadsell,
  M.~Balcan, and H.~Lin, Eds., vol.~33.\hskip 1em plus 0.5em minus 0.4em\relax
  Curran Associates, Inc., 2020, pp. 3094--3105.

\bibitem{hastie2017elements}
T.~Hastie, R.~Tibshirani, and J.~Friedman, ``The elements of statistical
  learning: data mining, inference, and prediction,'' 2017.

\bibitem{liang2022operation}
Z.~Liang, R.~Mieth, and Y.~Dvorkin, ``Operation-adversarial scenario
  generation,'' \emph{Electric Power Systems Research}, vol. 212, p. 108451,
  2022.

\bibitem{Park_Van_Hentenryck_2023}
S.~Park and P.~Van~Hentenryck, ``Self-supervised primal-dual learning for
  constrained optimization,'' \emph{Proceedings of the AAAI Conference on
  Artificial Intelligence}, vol.~37, no.~4, pp. 4052--4060, Jun. 2023.

\bibitem{pmlr-v63-Gao09}
W.~Gao, X.-Y. Niu, and Z.-H. Zhou, ``Learnability of non-i.i.d.'' in
  \emph{Proceedings of The 8th Asian Conference on Machine Learning}, ser.
  Proceedings of Machine Learning Research, R.~J. Durrant and K.-E. Kim, Eds.,
  vol.~63.\hskip 1em plus 0.5em minus 0.4em\relax The University of Waikato,
  Hamilton, New Zealand: PMLR, 16--18 Nov 2016, pp. 158--173.

\bibitem{6496182}
C.~Gong, X.~Wang, W.~Xu, and A.~Tajer, ``Distributed real-time energy
  scheduling in smart grid: Stochastic model and fast optimization,''
  \emph{IEEE Transactions on Smart Grid}, vol.~4, no.~3, pp. 1476--1489, 2013.

\bibitem{9585298}
S.~Gao, C.~Xiang, M.~Yu, K.~T. Tan, and T.~H. Lee, ``Online optimal power
  scheduling of a microgrid via imitation learning,'' \emph{IEEE Transactions
  on Smart Grid}, vol.~13, no.~2, pp. 861--876, 2022.

\bibitem{still2018lectures}
G.~Still, ``Lectures on parametric optimization: An introduction,''
  \emph{Optimization Online}, p.~2, 2018.

\bibitem{milgrom2002envelope}
P.~Milgrom and I.~Segal, ``Envelope theorems for arbitrary choice sets,''
  \emph{Econometrica}, vol.~70, no.~2, pp. 583--601, 2002.

\bibitem{mangasarian1987lipschitz}
O.~L. Mangasarian and T.-H. Shiau, ``Lipschitz continuity of solutions of
  linear inequalities, programs and complementarity problems,'' \emph{SIAM
  Journal on Control and Optimization}, vol.~25, no.~3, pp. 583--595, 1987.

\bibitem{ma2022lecture}
T.~Ma, ``Lecture notes for machine learning theory (cs229m/stats214),'' 2022.

\bibitem{kakade2008complexity}
S.~M. Kakade, K.~Sridharan, and A.~Tewari, ``On the complexity of linear
  prediction: Risk bounds, margin bounds, and regularization,'' \emph{Advances
  in neural information processing systems}, vol.~21, 2008.

\bibitem{bartlett2002rademacher}
P.~L. Bartlett and S.~Mendelson, ``Rademacher and gaussian complexities: Risk
  bounds and structural results,'' \emph{Journal of Machine Learning Research},
  vol.~3, no. Nov, pp. 463--482, 2002.

\bibitem{joukovsky2021generalization}
B.~Joukovsky, T.~Mukherjee, H.~Van~Luong, and N.~Deligiannis, ``Generalization
  error bounds for deep unfolding rnns,'' in \emph{Uncertainty in Artificial
  Intelligence}.\hskip 1em plus 0.5em minus 0.4em\relax PMLR, 2021, pp.
  1515--1524.

\bibitem{valleperez2020generalizationboundsdeeplearning}
\BIBentryALTinterwordspacing
G.~Valle-Perez and A.~A. Louis, ``Generalization bounds for deep learning,''
  2020. [Online]. Available: \url{https://arxiv.org/abs/2012.04115}
\BIBentrySTDinterwordspacing

\bibitem{truong2024rademachercomplexitybasedgeneralizationbounds}
\BIBentryALTinterwordspacing
L.~V. Truong, ``On rademacher complexity-based generalization bounds for deep
  learning,'' 2024. [Online]. Available: \url{https://arxiv.org/abs/2208.04284}
\BIBentrySTDinterwordspacing

\bibitem{9862539}
Q.~Hou, G.~Chen, N.~Dai, and H.~Zhang, ``Distributionally robust
  chance-constrained optimization for soft open points operation in active
  distribution networks,'' \emph{CSEE Journal of Power and Energy Systems}, pp.
  1--11, 2022.

\bibitem{10508977}
Q.~Hou, N.~Dai, and Y.~Huang, ``Voltage regulation enhanced hierarchical
  coordinated volt/var and volt/watt control for active distribution networks
  with soft open points,'' \emph{IEEE Transactions on Sustainable Energy},
  vol.~15, no.~3, pp. 2021--2037, 2024.

\end{thebibliography}

\begin{IEEEbiography}[{\includegraphics[width=1in,height=1.25in,clip,keepaspectratio]{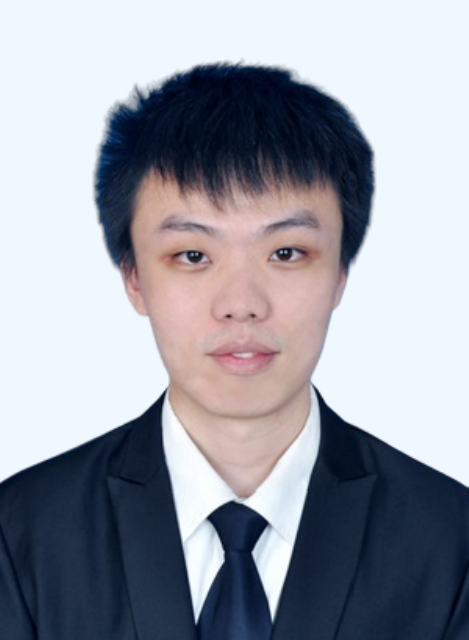}}]{Ge Chen}(S'20-M'23)
 received the B.S. degree from Huazhong University of Science and Technology, Wuhan, China, in 2015, and the M.S. degree from Xi'an Jiaotong University, Xi'an, China, in 2018, both in Thermodynamic Engineering. He earned his Ph.D. degree in Electrical and Computer Engineering from the University of Macau, Macao, China, in 2023. From 2023 to 2025, he was a Postdoctoral Research Associate at Purdue University. He is currently an Assistant Professor in the School of Engineering at Great Bay University, Dongguan, China. His research interests include include power system optimization and deep learning applications for intelligent decision-making.
\end{IEEEbiography}

\begin{IEEEbiography}[{\includegraphics[trim=1.8cm 0cm 1cm 0cm,width=1in,height=1.25in,keepaspectratio,clip]{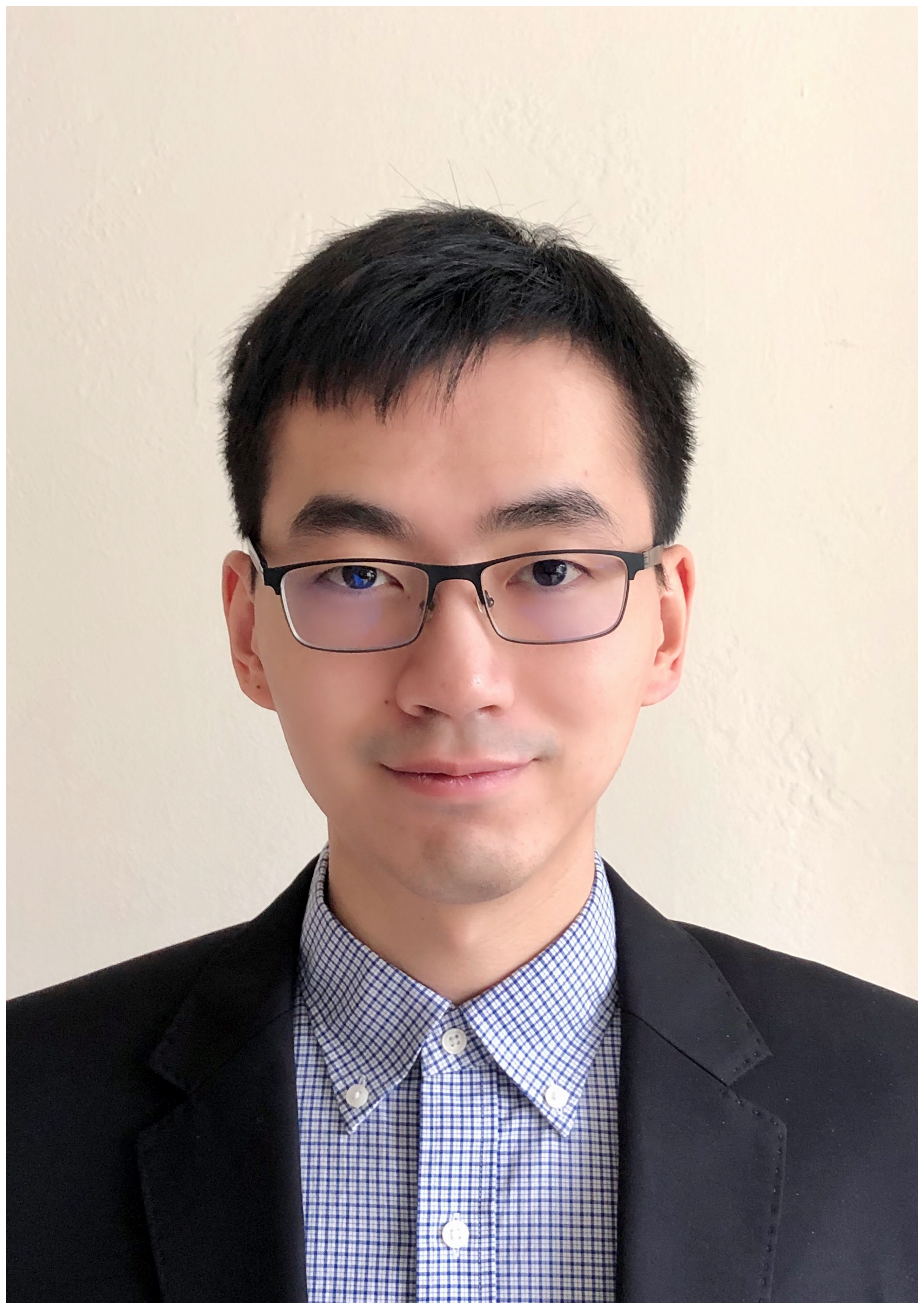}}]{Junjie Qin}
 is an Assistant Professor in the School of Electrical and Computer Engineering at Purdue University. Before joining Purdue, he was a postdoctoral researcher at University of California, Berkeley. He received a Bachelor of Engineering degree in Hydraulic and Hydropower Engineering (2010) and a Bachelor of Economics degree (2010) from Tsinghua University, Beijing, China. He obtained a Ph.D. degree in Computational and Mathematical Engineering (2018) from Stanford University, where he also received an M.S. degree in Civil and Environmental Engineering (2011) and an M.S. degree in Statistics (2017). His research interests include electric energy systems and transportation networks. He has received the NSF CAREER Award, the Google Research Scholar Award, the Satre Family Fellowship on Energy and Sustainability, and several best paper awards including the Energy System TC Outstanding Student Paper Award (as the advisor) from the IEEE Control and Systems Society in 2022, the O. Hugo Schuck Best Paper Award by the American Automatic Control Council in 2020, the Best Student Paper Award at the 23rd IEEE International Conference on Intelligent Transportation Systems in 2020, and the Best Student Paper Finalist at the 55th IEEE Conference on Decision and Control in 2016. 
 \end{IEEEbiography}

\end{document}